\documentclass[journal]{IEEEtran}
\usepackage{cite}
\usepackage[pdftex]{graphicx}
\usepackage{amsmath,amsthm,amssymb}
\usepackage{algorithmic}
\usepackage{algorithm}
\usepackage{color,soul}
\usepackage{eurosym}
\newtheorem{prop}{Proposition}
\usepackage{empheq}
\usepackage{tabularx}
\usepackage{float}
\usepackage{subfig}
\usepackage{booktabs}
\newcolumntype{C}[1]{>{\centering\arraybackslash}p{#1}}
\begin{document}
\title{Design of aggregators for the day-ahead management of microgrids providing active and reactive power services}

\author{Alessio La Bella,
        Marcello Farina,
        Carlo Sandroni,
        Riccardo Scattolini%
\thanks{Alessio La Bella, Marcello Farina and Riccardo Scattolini are with the Department of Electronics, Information and Bioengineering,  Politecnico di Milano, Via G. Ponzio 34/5, 20133 Milano (Italy). Email to:        {\tt\small name.surname@polimi.it}}
\thanks{Carlo Sandroni is with the Power Generation Technologies and Materials Department, Ricerca Sistema Energetico, Via R. Rubattino 54, 20134 Milano (Italy). Email to: {\tt\small carlo.sandroni@rse-web.it}}}
\maketitle

\begin{abstract}
\textit{The increasing diffusion of distributed energy generation systems requires the development of new control paradigms for the coordination of micro-generators, storage systems, and loads aimed at maintaining the efficiency and the safe operability of the electricity network. MicroGrids (MGs) are an interesting way to locally manage groups of generation devices, but they cannot singularly provide a significant contribution to sustain the main electricity grid in terms of ancillary services, such as the availability of a minimum amount of power reserve for the frequency regulation. For these reasons, in this paper we propose a framework for the aggregation and coordination of interconnected MGs to provide ancillary services to the main utility. The proposed framework is structured in three main phases. In the first one, a distributed optimization algorithm computes the day-ahead profile of the active power production of the MGs based on the available forecasts of the renewable sources production and the loads absorption. In this phase, scalability of the optimization problem and confidentiality requirements are guaranteed. In the second phase, reactive power flows are scheduled and it is ensured that the active power trends planned in the first phase do not compromise the voltage/current limitations. A final third phase is used to schedule the active and reactive power profiles of the generation units of each MG to make them consistent with the requirements and results of the previous two phases. The developed method is used for control of the IEEE 13-bus system network and the results achieved are thoroughly discussed in terms of performance and scalability properties.}
\end{abstract}


\section{Introduction}
\begin{figure}[t!]
\centering
\includegraphics[width=1\linewidth]{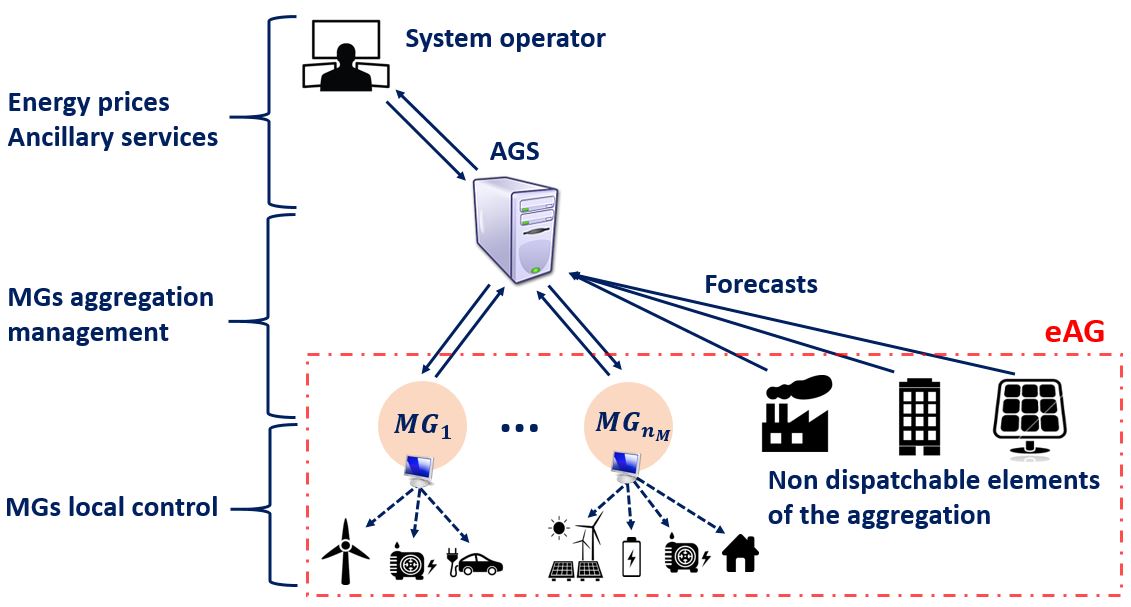}
\caption{Electrical aggregation management framework}
\label{fig:schema}
\end{figure}

\IEEEPARstart{T}HE wide diffusion of renewables, distributed energy generators, energy storage systems, and flexible loads has many advantages in terms of efficient and flexible management of the electricity grid and pollutants reduction. However, this evolution also poses new coordination problems which must be solved to guarantee reliable electrical services, like power production and frequency/voltage regulation. A partial and popular solution consists in grouping  generators, storages and loads in active small-scale grids, named MicroGrids (MGs), which can act as flexible energy consumers or providers, able to operate both when connected to the main grid and in isolated mode \cite{lasseter2002microgrids, guerrero2011hierarchical, la2017hierarchical}. In grid-connected mode, the internal generation units are usually managed in order to maximize the MG profit coming from the external energy trade, see \cite{parisio2014model, cominesi2016multi}, but, as the number of connected MGs increases, this selfish approach may cause serious problems at the grid level and soon also MGs, as well as smart buildings, will be required to provide external ancillary services, see  \cite{majzoobi2017application, joos2000potential, qureshi2016economic}. One of the main services that can be offered to the utility grid is the provision of active power reserve that can be later used to compensate frequency deviations \cite{yuen2007feasibility}.
This service, generally requested to traditional generation plants, is characterized by minimum requirements on the amount of power reserve to provide, and it can be hardly required to single MGs. Indeed, a single MG would have a negligible impact on the overall grid system and, moreover, it is commonly characterized by a limited capability, having also to satisfy the internal loads\cite{koliou2014demand,yuen2011provision}.
 A possible solution is to manage the available MGs in a collective fashion by aggregating interconnected MGs and non-dispatchable elements, and considering them as an unique electrical aggregation (eAG) controlled by an external entity, named aggregator supervisor (AGS), potentially capable to provide external ancillary services.
Considering eAGs allows both to reach an adequate size and to simplify the managing tasks in charge of the system operators, namely the Transmission System Operator (TSO) and the Distribution Network Operator (DNO), which can interact directly with the AGS system. The corresponding overall control structure is schematically represented in Figure \ref{fig:schema}. The aggregator supervisor can be designed according to centralized approaches, e.g. \cite{yuen2011provision, vrettos2016robust}, but these could be inapplicable in view of the large size of the overall system, leading to communication and scalability issues, and of confidentiality reasons, i.e., of possible restrictions to the information to be transmitted by the MGs, such as their internal costs,  generators characteristics and load demand profiles. Considering that MGs are locally managed by their own MGs Central Controllers (MGCC), the overall eAG scheduling process can be solved in a distributed manner, coordinating the MGs operations without requiring internal sensitive information. \\
In this paper we propose a distributed method for the design of an AGS planning the active power production/absorption of MGs on a day-ahead basis according to the daily prices of the sold/bought energy and to forecasts of loads and renewable energy sources. In this approach, the MGs autonomously manage their power sources while the AGS guarantees the best economic management of the eAG and that a minimum amount of power reserve is always allocated. Specifically, two kinds of active power reserve are considered: up and down power reserve capacity. They correspond to the power margins for increasing and decreasing the output power of the generation units in order to compensate severe frequency deviations. \\ An additional external service can be provided by using MGs as reactive power producers/consumers. This allows to regulate the line currents and voltages inside the eAG and to minimize the eAG power losses. This second task is addressed modelling each MG as a single equivalent generator. This allows to solve the reactive power scheduling problem in a centralized way and, at the same time, to preserve the MGs internal sensitive data. This is made possible since the reactive power service does not affect internal production costs as it involves just the regulation of the inverters and of the excitation generation systems. In other words, it can be centrally asked each MG to produce/absorb a certain amount of reactive power without affecting its internal economy. According to these considerations, the proposed scheduling method is structured as follows.\\ \\
\textit{Phase 1:} \textbf{Distributed active power and reserve dispatch.} \\
The day-ahead active power production/absorption profiles of the MGs generation units, minimizing the internal production costs and globally providing the minimum required active power reserve, are computed with a distributed optimization algorithm, the Alternating Direction Method of Multipliers (ADMM), see \cite{boyd2011distributed}. \\	
\\
\textit{Phase 2:} \textbf{Power flow feasibility and reactive power planning.}\\
Each MG is requested to communicate its overall reactive power capability to the AGS. The MGs' reactive power production/absorption is then scheduled to satisfy the voltage/current limits and to minimize the overall power losses inside the eAG. In this step, if an On Load Tap Changer transformer (OLTC) is available at the interconnection with the main utility, its voltage set-points can be also managed. In case the power flow feasibility is compromised, the AGS can also ask the MGs to vary the active power profiles scheduled in \textit{Phase 1} considering the line limitations and the network topology.
\\ \textit{Phase 3:}  \textbf{Final scheduling of the MGs generation units.}\\
In this final phase, each MG must schedule the active and reactive power daily trends of its generation units considering the requirements from \textit{Phase 1} and \textit{Phase 2}.\\
Two situations may occur:
\begin{itemize}
\item  \textit{Phase 3(a)}: the active power profiles scheduled in \textit{Phase 1} do not compromise the electrical feasibility, meaning that active power variations are not requested by \textit{Phase 2}. In this case, the MGs independently schedule the reactive power daily references for their local generators.
\item  \textit{Phase 3(b)}: Active power variations are needed to cope with current and voltage limitations. In this case, the MGs must not only schedule the reactive power daily production of their generators, but also re-schedule their active power production considering the overall set-points provided by the previous phases.
\end{itemize}
\begin{table}[t!]
\caption{Acronyms}
\label{tab:acro}
\centering
\begin{tabular}{l l}
		\hline
	Acronym  &  Description\\ \hline \vspace{0.8mm}
     ADMM & alternating direction method of multipliers \\ \vspace{0.8mm}
     AGS & aggregator supervisor\\ \vspace{0.8mm}
     DNO & distribution network operator\\
     \vspace{0.8mm}
     TSO & transmission system operator\\ \vspace{0.8mm}
	 eAG &  eletrical aggregation \\ \vspace{0.8mm}
	 MG & microgrid\\ \vspace{0.8mm}
	 OLTC & on load tap changer  \\   \vspace{0.8mm}
\end{tabular}
\end{table}
%

 \begin{figure}[t!]
\centering
\includegraphics[width=1\linewidth]{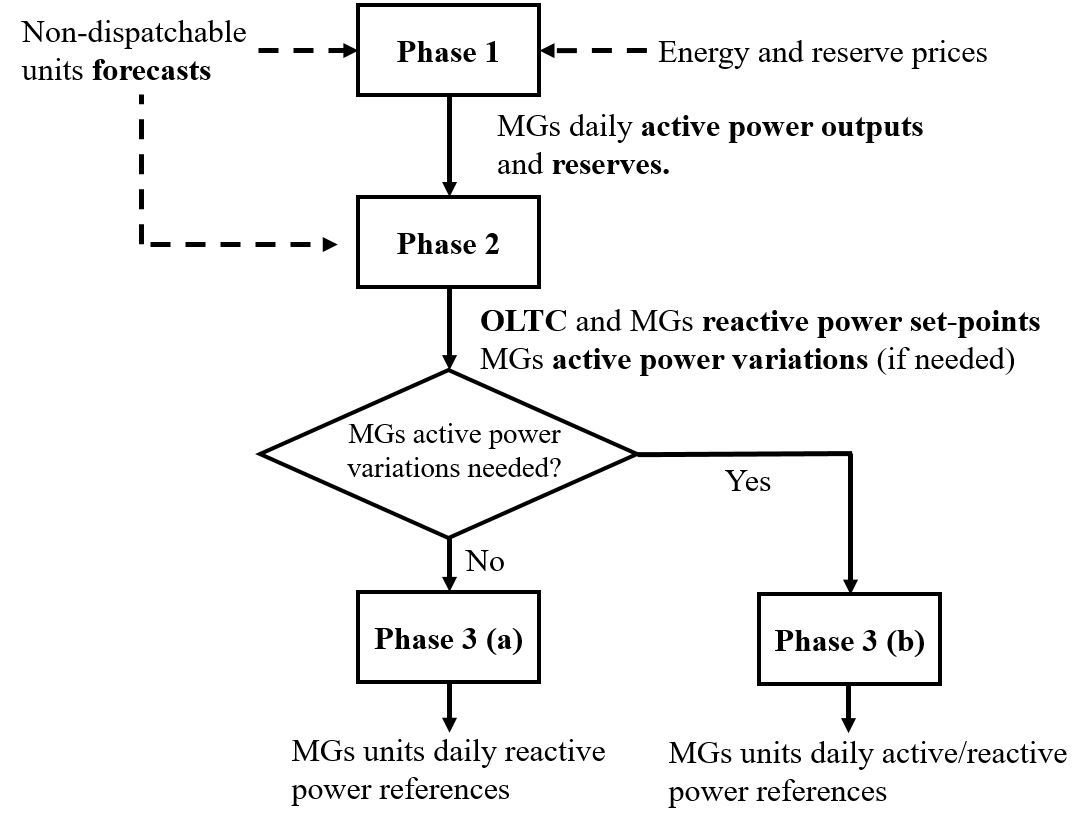}
\caption{Flowchart of the eAG day-ahead scheduling process}
\label{fig:flowchart}
\end{figure}
\bigskip
A flowchart of the algorithm describing the inputs and outputs of each phase is shown in Figure \ref{fig:flowchart}.
As already mentioned, this three-step procedure is supposed to be performed off-line, relying on the day-ahead energy prices and on the forecasts of the load and weather trends, supposed to be available from historical data. In this way, the AGS manages to schedule the eAG power output achieving the best internal economic management of the MGs respecting both the frequency reserve and the electrical requirements, and then it can communicate it directly to the system operators, e.g. TSO and DNO.

\pagebreak
The proposed approach is different from the ones discussed in \cite{balandat2014contract} and \cite{mai2015economic}, where the goal of the aggregator supervisor was to maximize its own profit, acting as an intermediate between the system operators and the MGs; moreover, in \cite{balandat2014contract} and \cite{mai2015economic} a centralized framework was adopted. In \cite{gkatzikis2013role,kim2011two,dagdougui2014decentralized, brehm2017distributed}, distributed optimization frameworks have been designed to manage aggregations of energy resources and loads, without considering the ancillary services provision. A distributed management of commercial buildings for frequency reserve provision has been considered in \cite{rey2018strengthening} where however, a simplified thermal modelling of buildings has been adopted. Differently from the mentioned contributions, in this paper the aggregator supervisor is also entitled of ensuring the electrical feasibility of the scheduled power flows, performing the Phase 2. 
Concerning this optimization step, a distributed approach, as considered in \cite{guo2017case, magnusson2015distributed}, can be critical due to the nonlinearity and non-convexity of the underlying power flow equations,  which negatively impact on the convergence of most distributed optimization-based algorithms \cite{molzahn2017survey}. However, as already stated, distributed implementations are not necessary, and a centralized approach can be used still preserving confidentiality requirements and non affecting internal production costs.\\
This paper extends the preliminary work reported in \cite{labellaECC2018}, where just the reserve provision was considered and a low-effective distributed decomposition method was used. The paper is structured as follows: in Section \ref{sec:2}, the MGs modelling is presented and the proposed distributed optimization algorithm used in Phase 1 is described. The reactive power management aimed at guaranteeing the feasibility of the proposed solution is presented in Section \ref{sec:3}. In Section \ref{sec:4}, the final optimization performed in Phase 3 is described. The  results achieved by applying the proposed framework to the IEEE 13-bus network, configured to include MGs, are extensively discussed in Section \ref{sec:5}. Finally, some conclusions are drawn in Section \ref{sec:6}.
\begin{table}[t!]
\caption{Optimization variables and system parameters}
\label{tab:MG_active_opt_var}
\centering
\begin{tabular}{l l}
		\hline 
	Symbol  &  Description\\ \hline \vspace{1mm}
		 \(p^{\,g},q^{\,g}\) & MG generator active/reactive power set-points [kW, kVar]\\ \vspace{1mm}
	 \(p^{\,b},q^{\,b}\) &  MG battery active/reactive power set-points [kW, kVar] \\ \vspace{1mm}
	 $p^{l},q^{l}$ & MG load active/reactive power forecast [kW,kVar]\\ \vspace{1mm}
$p^{r}$ & MG renewable source active power forecast [kW]\\ \vspace{1mm}
$q^{r}$ & MG renewable source reactive power set-point [kVar]\\ \vspace{1mm}
	    	$ p^{\scriptstyle{mg}}, q^{\scriptstyle{mg}}$ & MG active/reactive output power [kW, kVar]  \\ \vspace{1mm}
	\(s^{b}\) & MG battery state of charge (SOC) [\%] \\ \vspace{1mm}
$ C^{\,b} $ & MG battery capacity [kWh] \\ \vspace{1mm}
$ c^{\,b} $ & MG battery usage cost [\euro/kWh$^2$]\\ \vspace{1mm}
 ${a^{g}}, {b^{g}}, {c^{g}} $  & MG generator cost coefficients [\euro/kWh$^2$,\,\euro/kWh,\,\euro] \\  \vspace{1mm}
$ r_p^{\,g  \uparrow} ,r_p^{\,g \downarrow} $ & MG generator up/down active power reserves  [kW]\\ \vspace{1mm}
$ r_p^{\,b \uparrow},r_p^{\,b \downarrow}  $ & MG battery up/down active power reserves  [kW] \\ \vspace{1mm}
    $ {r_{p}^{\scriptstyle{mg}\uparrow}} ,  {r_{p}^{\scriptstyle{mg}\downarrow}}$ \hspace*{-2mm}& MG total up/down active power reserves  [kW]\\  \vspace{1mm}
     $ {r_{q}^{\scriptstyle{mg}\uparrow}} ,  {r_{q}^{\scriptstyle{mg}\downarrow}}$ \hspace*{-2mm}& MG total up/down reactive power reserves  [kVar]\\  \vspace{1mm}
    $\rho_e^s,\, \rho_e^b$ & Grid energy selling and buying price [\euro/kWh]	\\ \vspace{1mm}
$\rho_{r_p}^{\uparrow},\,\rho_{r_p}^{\downarrow}$ & Up/down capability  active reserve prices [\euro/kWh]	\\
   ${p}^{\textit{\tiny{L}}},{q}^{\textit{\tiny{L}}}$  & Active/reactive power forecast of an eAG\\  \vspace{1mm} & absorption node [kW,kVar] \\
${p}^{\textit{\tiny{R}}}$  & Active power forecast an eAG non-dispatchable\\  \vspace{1mm} & generation node [kW]	\\ \vspace{1mm}
  $ {p}^{\text{\tiny{AG}}}$ &  eAG total active output power [kW] \\ \vspace{1mm}
 $ {r}_p^{\text{\tiny{AG}}\uparrow}, {r}_p^{\text{\tiny{AG}}\downarrow}$&  eAG total active up/down active power reserves [kW]\\ \vspace{1mm}
\( V^o \)         & Nominal network voltage [V]\\ \vspace{1mm}
	    	$ \Delta V^s $ & OLTC voltage step-size with respect to $V^o$ [\%]\\ \vspace{1mm}
     $ \Delta \bar{V}^s $ & Maximum OLTC voltage variation  [\%] \\ \vspace{1mm}
		 		\( \theta^s \)         & OLTC step selector [boolean]\\ \vspace{1mm}
      \(V  \) & Nodal voltage magnitude  [V]\\ \vspace{1mm}
$P^{\,loss}$ & Active power line loss [kW] \\ \vspace{1mm}
$I$ & Line current magnitude [A] \\
\hline
\end{tabular}
\end{table}
\section{Phase 1: active power dispatch and reserve provision}\label{sec:2}
The network interconnecting the eAG elements can be modelled as a bi-directional graph $\mathcal{G(N,E)}$ with nodes $\mathcal{N}= \{1,\hdots,n\}$ and edges $\mathcal{E} \subseteq \mathcal{N} \times \mathcal{N}$. The eAG includes a set $\mathcal{N_M}$ of flexible nodes (e.g. the MGs), a set $\mathcal{N_L}$ of non-dispatchable load absorption nodes (e.g., commercial buildings), and a set $\mathcal{N_R}$ of non-dispatchable generation nodes (e.g., renewable sources power plants), where $\mathcal{N_M}$, $\mathcal{N_L}$, $\mathcal{N_R} \subseteq  \mathcal{N}$. All the variables and parameters that will be used in the next sections are described in Table \ref{tab:MG_active_opt_var}. As a convention, all the power values are defined to be positive if delivered and negative if absorbed. Moreover, the maximum and minimum limits of each variable are denoted with a bar over or below the variable, respectively.
\subsection{Microgrid modelling and problem formulation}\label{subsec:mg_model}
In \emph{Phase 1} the economic management of the internal resources is performed by computing, in a distributed way, the optimal day-ahead profiles of the active power production of the MGs. For this purpose MGs are modelled as discrete time systems with sampling time $\tau=15$ min, being grid energy prices and weather forecasts usually provided with the same time rate. In view of the adopted time horizon of 24 hours, $N=96$ steps must be considered in the optimization problem.\\
The $i$-th microgrid, denoted as   MG$_i$ with $i \in \mathcal{N_M}$, is generally equipped with $n_i^g$ fuel-based generators, $n_i^b$ storage units, $n_i^r$ renewable generators and has a given load demand. Here, it is assumed that the loads are neither shiftable nor interruptible.
In the timescale defined by the sampling time $\tau$, the dynamics of the generation units can be neglected. The active power capability constraints are:
\begin{align}
\underline{p}^{\,g}_{j_i} \, \leq \, &{p}^{\,g}_{j_i}(t) \, \leq \,  \bar{p}^{\,g}_{j_i} \label{gen_limits}  \\
\underline{p}^{\,b}_{k_i}  \leq \, &{p}^{\,b}_{k_i}(t) \, \leq \,  \bar{p}^{\,b}_{k_i} \label{batt_limits}
\end{align}
where $t$ is the time index, while $j_i\in (1 ,\hdots,n^g_i)$ and $ k_i \in (1 ,\hdots,n^b_i)$ represent the $j$-th generator and the $k$-th battery installed in MG$_i$, respectively. The dynamics of the state of charge (SOC) of the batteries is:
\begin{align}
s^{b}_{k_i}(t+1) \; = \; &s^{b}_{k_i}(t) \; - \; 100\,\frac{\tau}{C^{\,b}_{k_i} } \, p^{b}_{k_i}(t)
\label{SOC_dyn}
\end{align}
where for simplicity, the charge/discharge efficiencies have been neglected.
The SOC must be constrained between minimum and maximum bounds; moreover, it is often required that, for each battery, the amount of energy stored at the end of the day equals the one at the beginning:
\begin{align}
\begin{split}
\underline{s}^{\,b}_{\,k_i} \, \leq \; {s}^{\,b}_{\,k_i}&(t) \, \leq \,  \bar{s}^{\,b}_{\,k_i}   \\
{s}^{\,b}_{\,k_i}(N) = &\;{s}^{\,b}_{\,k_i}(0)
\end{split}
\label{SOC_constr}
\end{align}
The overall MG output active power is given by the following internal power balance; it is expected that the output active power of  MG$_i$ verifies some bounds, dictated for instance by an electrical contract with the DNO:
\begin{flalign}
 p_i^{\scriptstyle{mg}}(t)\;=  \sum_{k_i=1}^{n^b_{i}}p^{\,b}_{k_i}(t)& + \sum_{j_i=1}^{n^g_{i}}p^{\,g}_{j_i}(t)  + \sum_{p_i=1}^{n^r_{i}}p^{\,r}_{p_i} + p^{l}_{i}(t) \label{mg_output_power}&&\\
\underline{p}_i^{\scriptstyle{mg}} &\leq   p_i^{\scriptstyle{mg}}(t)\leq \bar{p}_i^{\scriptstyle{mg}}&&
\label{mg_output_constr}
\end{flalign}
The active power reserves provided by the fuel-based generators correspond to the remaining power margins with respect to the capability limits:
\begin{align}
r^{\,g \uparrow}_{p,j_i} (t) &=  \bar{\,p}^{g}_{\,j_i} \; - \; {p}^{\,g}_{\,j_i}(t) \\
\label{gen_reserve}
r^{\,g \downarrow}_{p,j_i} (t) &=   \, {p}^{g}_{\,j_i}(t)  \; - \; \underline{p}^{g}_{\,j_i}
\end{align}
As for the active power reserves of the batteries, not only the capability limits must be considered, but also the amount of stored energy. Indeed, the battery active power bounds are defined as the minimum values between the capability limits and the maximum amount of generating/absorbing power based on the actual SOC:
\begin{align*}
r^{b \, \uparrow}_{p,k_i} (t)  &=  \text{min} \, \{ \;  \; \;\;\,\bar{p}^b_{k_i}, \;\;  \frac{ \big( s^b_{k_i}(t) \, - \, \underline{s}^b_{k_i}\,\big)}{100} \frac{ C^b_{k_i}}{ \tau} \;\} \; - \; {p}^b_{k_i}(t)  \\
r^{b \, \downarrow}_{p,k_i} (t)  &=  \text{min} \, \{ \; - \; \underline{p}^{\,b}_{k_i}\, , \; \; \frac{\big( \, \bar{s}^b_{k_i} \, - \, {s}^b_{k_i}(t) \,\big)}{100} \frac{ C^b_{k_i}}{ \tau} \}\; + \; {p}^{\,b}_{k_i}(t)
\end{align*}
The nonlinearity of these expressions negatively affects the complexity of the optimization problem. Therefore, they are reformulated  as follows:
\begin{align}
r^{b \uparrow}_{p,k_i} (t)   + \; {p}^b_{k_i}(t)\,  &\leq \; \;\bar{p}^b_{k_i}
\label{batt_upreserve_costr_1} \\
\quad \; \, r^{b\uparrow}_{p,k_i} (t)  + \; {p}^b_{k_i}(t)\, &\leq\; \frac{ \big( s^b_{k_i}(t) \, - \, \underline{s}^b_{k_i}\,\big)}{100} \frac{ C^b_{k_i}}{ \tau}
\label{batt_upreserve_costr_2} \\
 r^{b \downarrow}_{p,k_i} (t) -\, {p}^b_{k_i}(t)  \,  &\leq\;   - \, \underline{p}^b_{k_i}
\label{batt_downreserve_costr_1} \\
 r^{b\downarrow}_{p,k_i} (t)  - \,{p}^b_{k_i}(t)  \, &\leq \;\frac{\big( \, \bar{s}^b_{k_i} \, - \, {s}^b_{k_i}(t) \,\big)}{100} \frac{ C^b_{k_i}}{ \tau}
\label{batt_downreserve_costr_2}
\end{align}
As it will be shown, the amount of active power reserve is a gain for a MG. This implies that the optimal values of the reserves of the batteries are the minimum between the
admissible values in expression \eqref{batt_upreserve_costr_1}-\eqref{batt_upreserve_costr_2} for the up reserve and  \eqref{batt_downreserve_costr_1}-\eqref{batt_downreserve_costr_2} for the down reserve.
The overall reserves are obtained by summing all the contributions of each MG. It is assumed that power production of the renewable generators can be used as down active power reserve since
it can be reduced in case of frequency unbalances. Formally, we define
\begin{align}
 {r_{p}^{\scriptstyle{mg}\uparrow}}   (t&) =  \sum_{k_i=1}^{n^b_{i}}r^{\,b \uparrow}_{p,k_i}(t) \; + \;\sum_{j_i=1}^{n^g_{i}}r^{\,g\uparrow}_{p,j_i}(t) \label{mg_up_reserve}\\
  {r_{p}^{\scriptstyle{mg}\downarrow}}  (t&) =  \sum_{k_i=1}^{n^b_{i}}r^{\,b \downarrow}_{p,k_i}(t)  + \;\sum_{j_i=1}^{n^g_{i}}r^{\,g\downarrow}_{p,j_i}(t) + \sum_{p_i=1}^{n^b_{i}}p^{r}_{p_i}(t) \label{mg_dwn_reserve}
\end{align}
The goal of the MG internal management is to minimize the production costs and maximize the gain from the external trading. Therefore, the cost function for the generic MG$_i$ management is:
\begin{align}
\begin{split}
J^{mg}_i =  \,  & \sum_{t=1}^{N}\,\underbrace{\sum_{j_i=1}^{n^g_{i}}({a}_{j_i}^{\,g}\tau^2({p^{\,g}_{j_i}}(t))^2  + b_{j_i}^{\,g}\,\tau\, p^{\,g}_{j_i}(t)+  {c}_{j_i}^{\,g})}_{\alpha}  \,+ \\ + & \sum_{t=2}^{N} \underbrace{\sum_{k_i=1}^{n^b_{i}}c^{\,b}_{k_i} \tau^2\,( p^{b}_{k_i}(t)  - p^{b}_{k_i}(t-1) )^2}_{\beta} \, +\\ - & \sum_{t=1}^{N}\,\underbrace{h^p_i\,(\, p_i^{\scriptstyle{mg}}(t), \,  {r_{p}^{\scriptstyle{mg}\uparrow}}(t), \,  {r_{p}^{\scriptstyle{mg}\downarrow}} (t) \, )}_{\eta}
\end{split}
\label{muG_costfunction}
\end{align}

In \eqref{muG_costfunction}, $\alpha$ represents the  MG generators fuel cost, expressed as a second-order polynomial function with respect to the generated power \cite{sonmez2013estimation}. The term $\beta$ includes the squared variation of the active power set-point of the batteries to avoid their frequent and excessive charge and discharge (note that that the costs are multiplied by $\tau$ since they are energy costs). Finally, the function $h^p_i$ in $\eta$ represents the gain/cost of the MG considering the external trade, which will be defined later since it depends on the AGS management.
MG$_i$ is compactly represented by defining the following vectors: \\[3mm] ${x^{p}_i}=[( p^g_{j_i},r^{g\uparrow}_{p,j_i}, r^{g\downarrow}_{p,j_i})_{\forall j_i \in \{1,n_i^g\}}, ( p^b_{k_i} ,r^{b\uparrow}_{p,k_i}, r^{b\downarrow}_{p,k_i})_{\forall k_i \in \{1,n_i^b\}}]'$,\\\\
${d^p_{i}}=[ \,p^{l}_i \,,\, p^r_{\forall p_i \in  \{1,n^r_i\}}\,]'$,  \\\\
${y^p_{i}}=[ \,p_i^{\scriptstyle{mg}}, \,  {r_{p}^{\scriptstyle{mg}\uparrow}}, \,  {r_{p}^{\scriptstyle{mg}\downarrow}} \,]'$, \\[3mm] where $x^p_{i}$ includes the MG$_i$ internal active power variables, $d^p_{i}$  is the vector of non-dispatchable active power trends and $y^p_{i}$ is the vector of the MG$_i$ output variables.  For compactness, the variables referred to the whole time horizon are expressed with bold symbols, e.g. $\mathbf{x}_i^\mathbf{p}=[x^p_{i}(1)',\hdots,x^p_{i}(N)']'$. Similarly, in compact form the MG$_i$ cost function can be rewritten as
\begin{subequations}
\begin{align}
J_i^{mg}={f}_i^{p}(\mathbf{x}_i^\mathbf{p}) \, - \,  {h}_i^{p}(\mathbf{y}_i^\mathbf{p}) 
\label{muG_costfunction_compact}
\end{align}
where ${f}_i^{p}$ includes the production cost of  MG$_i$, i.e. the terms $\alpha$ and $\beta$ in \eqref{muG_costfunction}. The MG$_i$ constraints  are compactly expressed as
\begin{align}
 \boldsymbol{A}_i^\mathbf{p}\,
\mathbf{x}_i^\mathbf{p} \, \leq\, \boldsymbol{b}_i^\mathbf{p} \label{mg_constr}
\end{align}
where $\boldsymbol{A}_i^\mathbf{p}$ and $\boldsymbol{b}_i^\mathbf{p}$ are properly defined. Finally, the outputs of MG$_i$ can be written as   
\begin{align}
&\,\mathbf{y}_i^\mathbf{p}\, = \,  \boldsymbol{C}_i^\mathbf{p} \,\mathbf{x}_i^\mathbf{p} \, + \,  \boldsymbol{M}_i^\mathbf{p}\, \mathbf{d}_i^\mathbf{p} \quad \label{mg_output}
\end{align}
with a proper definition of the matrices $\boldsymbol{C}_i^\mathbf{p}$ and $\boldsymbol{M}_i^\mathbf{p}$.
\label{MG_opt_problem}
\end{subequations}
\\
For the description of the AGS scheduling process, the MG$_i$ compact optimization problem  \eqref{MG_opt_problem} will be considered. This is not only related to the notational compactness, but also because MGs can be characterized by different cost functions, constraints and composition of generations units and this must not have an impact on the AGS management.
It is also worth noting that the ancillary services requirements, e.g. the provision of a minimum amount of power reserve, are not included in the  optimization problem \eqref{MG_opt_problem}. These constraints will be defined at the AGS management level.
\subsection{AGS problem formulation and distributed algorithm}
At a collective level, the AGS must consider the selling/buying energy prices and the gain for the provided active power reserve. The centralized AGS cost function is defined as follows:
\begin{subequations}
\begin{align}
\begin{split}
\mathbf{J^{\textbf{\tiny{AG}}}} &= \sum_{\forall i\in \mathcal{N_M}} f_i^p(\mathbf{x}_i^{\mathbf{p}}) \;- \;\tau \boldsymbol{\rho_e^{s}}{'}\, {\max}\,(\mathbf{{p}^{\textbf{\tiny{AG}}}} ,\boldsymbol{0}) \;+\\&+\;\tau\boldsymbol{\rho_e^{b}}{'}\, {\max}(- \mathbf{{p}^{\textbf{\tiny{AG}}}},\boldsymbol{0}) -\tau\boldsymbol{\rho_{r_p}^{\uparrow}}\hspace{-1mm}{'}  \, \mathbf{{r}_p^{\textbf{\tiny{AG}}\uparrow}} - \tau \boldsymbol{\rho_{r_p}^{\downarrow}}\hspace{-1mm}{'}  \, \mathbf{{r}_p^{\textbf{\tiny{AG}}\downarrow}}
\end{split}\label{AGS_centr_cf}
\end{align}
Note that the eAG at time $t$ can either absorb (${p}^{\text{\tiny{AG}}}(t)<0$) or release (${p}^{\text{\tiny{AG}}}(t)>0$) power. The functions ${\max}\,(\pm \mathbf{{p}^{\text{\tiny{AG}}}},\mathbf{0})$ in \eqref{AGS_centr_cf} compute the element-wise maximum between the eAG output power profile $\pm\mathbf{{p}^{\text{\tiny{AG}}}}$ and the zero column vector, in order to differently price the output eAG power considering either the selling or the buying energy prices. Moreover, the overall cost function considers also the gain for the provided up and down active power reserves by the whole eAG. \\
The constraints \eqref{mg_constr}-\eqref{mg_output} must be considered for all $i \in \mathcal{N_M}$, together with constraints for the minimum required up and down power reserve 
\begin{align}
\mathbf{{r}_p^{\textbf{\tiny{AG}}\uparrow}} \geq \; \mathbf{\underline{r}_p^{\textbf{\tiny{AG}}\uparrow}}\, ,\;\;  \mathbf{{r}_p^{\textbf{\tiny{AG}}\downarrow}} \geq \; \mathbf{\underline{r}_p^{\textbf{\tiny{AG}}\downarrow}}
 \label{AG_compact_centr_reserve_min}
\end{align}
Finally, the following expressions are stated
\begin{align}
&\mathbf{{r}_p^{\textbf{\tiny{AG}}\uparrow}}  = \sum_{\forall i \in \mathcal{N_M}}  \mathbf{r}_{\mathbf{p},\,i}^{\mathbf{\scriptscriptstyle{mg}\uparrow}},\quad\mathbf{{r}_p^{\textbf{\tiny{AG}}\downarrow}}  = \hspace{-2mm} \sum_{\forall i \in \mathcal{N_M}}\hspace{-2mm} \mathbf{r}_{\mathbf{p},\,i}^{\mathbf{\scriptscriptstyle{mg}\downarrow}} + \hspace{-2mm}\sum_{\forall j_r \in \mathcal{N_R}}  \hspace{-2mm}\mathbf{p}^{\textbf{\tiny{R}}}_{j_r} \label{AG_compact_centr_reserve} \\
&\mathbf{{p}^{\textbf{\tiny{AG}}}}  = \sum_{\forall i \in \mathcal{N_M}}\mathbf{p}_{i}^{\mathbf{\scriptscriptstyle{mg}}}\;+\; \sum_{\forall j_l \in \mathcal{N_L}} \mathbf{p}^{\textbf{\tiny{L}}}_{j_l} \;+ \;\sum_{\forall j_r \in \mathcal{N_R}} \mathbf{p}^{\textbf{\tiny{R}}}_{j_r}
\label{AG_compact_centr_powerbalance}
\end{align}
\label{AG_compact}
\end{subequations}
defining the total up and down eAG power reserve and the active power balance for the whole eAG.
From \eqref{AG_compact}, it is clear
that the AGS objective is to minimize the MGs production costs and maximize the overall profit due to the external trade 
with the main grid respecting the MGs internal contraints and the minimum up and down eAG power reserve constraints.
The overall optimization problem can be stated in the following compact form:
\begin{subequations}
\begin{align}
&\min_{\substack{ \,\mathbf{x}_{\forall i}^\mathbf{p},\,\mathbf{y}_{\forall i}^\mathbf{p} \\ \mathbf{z}^\mathbf{p}  }}\{\, \sum_{\forall  i\in \mathcal{N_M}} f_i^p(\mathbf{x}_i^\mathbf{p}) \; + \;g_z^p(\mathbf{z}^\mathbf{p} ) \;\}
\label{AG_cmpct_globalcf}\\
&\hspace{-15mm} \text{subject to} \nonumber\\
\begin{split}
& \boldsymbol{A}_i^\mathbf{p}\,
\mathbf{x}_i^\mathbf{p} \, \leq\, \boldsymbol{b}_i^\mathbf{p}\\
&\,\mathbf{y}_i^\mathbf{p}\, = \,  \boldsymbol{C}_i^\mathbf{p} \,\mathbf{x}_i^\mathbf{p} \, + \,  \boldsymbol{M}_i^\mathbf{p}\, \mathbf{d}_i^\mathbf{p}
\end{split}
 \qquad \forall i \in \mathcal{N_M}
\label{AG_cmpct_mgconstr}
\\
 &\, \mathbf{z}^\mathbf{p}    \; \in \mathcal{Z}^p \label{AG_cmpct_AGconstr}
 \\
 &\, \mathbf{z}^\mathbf{p}  \;- \sum_{\forall i\in \mathcal{N_M}} \mathbf{y}_i^\mathbf{p}  \;  - \,\mathbf{d}^\mathbf{p}  \;= \;\mathbf{0}\label{AG_cmpct_coupling}
\end{align}
\label{AG_cmpct}
\end{subequations}
where $\mathbf{z}^\mathbf{p} =[\mathbf{p}^{\textbf{\tiny{AG}}}\,{'},\mathbf{{r}_p^{\textbf{\tiny{AG}}\uparrow}}\,{'},\mathbf{{r}_p^{\textbf{\tiny{AG}}\downarrow}}\,{'}]\,'$ and
\begin{align*}
g_z^p(\mathbf{z}^\mathbf{p}) =&- \;\tau \boldsymbol{\rho_e^{s}}{'}\, {\max}\,(\mathbf{{p}^{\textbf{\tiny{AG}}}} ,\boldsymbol{0}) \;+ \;\tau\boldsymbol{\rho_e^{b}}{'}\, {\max}(- \mathbf{{p}^{\textbf{\tiny{AG}}}},\boldsymbol{0})\;+\\ & -\;\tau\boldsymbol{\rho_{r_p}^{\uparrow}}\hspace{-1mm}{'}  \, \mathbf{{r}_p^{\textbf{\tiny{AG}}\uparrow}}\; -\; \tau \boldsymbol{\rho_{r_p}^{\downarrow}}\hspace{-1mm}{'}  \, \mathbf{{r}_p^{\textbf{\tiny{AG}}\downarrow}}
\end{align*}
The eAG reserve requirements \eqref{AG_compact_centr_reserve_min} have been condensed in  \eqref{AG_cmpct_AGconstr} properly defining the set $\mathcal{Z}^p \subset \mathbb{R}^{3N,1}$, while the constraints \eqref{AG_compact_centr_reserve}-\eqref{AG_compact_centr_powerbalance} have been compacted in \eqref{AG_cmpct_coupling}, introducing the vector $\mathbf{d}^\mathbf{p}$  combining the eAG non-dispatchable output powers.
The formulation \eqref{AG_cmpct} clearly shows that the optimization problem includes the coupling constraints \eqref{AG_cmpct_coupling} which collects the MGs optimization variables and the AGS ones. Finally, we identify with $\mathbf{J^{\textbf{\tiny{AG}},*}}$ the optimal value of the cost function \eqref{AG_cmpct_globalcf}, supposing to  be able to solve \eqref{AG_cmpct} in a centralized fashion.	\\
However, as previously discussed, its centralized implementation is not advisable.
 As already mentioned, one of the main drawbacks of the optimization problem stated in \eqref{AG_cmpct} is that the AGS should know everything about the MGs internal structure, units characteristics and costs such that it can directly provide the set-points for their units. To overcome this and to confer scalability to the solution of the problem \eqref{AG_cmpct}, it is here proposed to rely on the Alternating Direction
Method of Multipliers (ADMM), a powerful distributed optimization algorithm with enhanced convergence properties with respect to the standard \textit{dual-decomposition} approaches \cite{boyd2011distributed}. To reach  convergence, ADMM requires the convexity of the cost function of the centralized problem. To this regard, considering \eqref{AG_cmpct_globalcf}, the terms $f_i(\mathbf{x}_i^\mathbf{p})$ have been defined to be convex, while, regarding $g_z^p(\mathbf{z}^\mathbf{p})$, its piece-wise structure does not directly guarantee any convexity properties. However, the following result can be stated (the proof is reported Appendix \ref{Prop1_proof}):
\begin{prop}
If $\rho_e^b(t) \geq \rho_e^s(t)$ $\; \forall t \, \in \, \{1, \hdots,N\}$, then the function $g_z^p(z^p(t))$ is convex (but not strictly).
\label{Prop1}
\end{prop}
The requirement that, at each time instant, the buying price must be greater than the selling one is realistic and also quite common in the literature, e.g. \cite{zhang2013robust,cominesi2018two}. Indeed, in case the eAG absorbs active power, the TSO/DNO afford some fixed costs to bring the energy to the eAG additionally to the generation ones, due for example to transmission or power losses. This is not the case for the eAG when the energy is produced \textit{in loco} by the MGs resources and sold to the main grid.
Therefore, the ADMM method can be applied, defining the following \textit{Augmented Lagrangian} function:
\begin{align}
\begin{split}
&L_{\mu}(\mathbf{x}_{\forall i}^\mathbf{p},\mathbf{y}_{\forall i}^\mathbf{p}, \mathbf{z}^\mathbf{p} ,\boldsymbol{\lambda}) =\hspace{-2mm} \sum_{\forall i\in \mathcal{N_M}}\hspace{-3mm} f_i^p(\mathbf{x}_i^\mathbf{p} )+g^p_z (\mathbf{z}^\mathbf{p} )   + \boldsymbol{\lambda}'( \mathbf{z}^\mathbf{p}  \,+\\&- \sum_{\forall i\in \mathcal{N_M}} \mathbf{y}_i^\mathbf{p}    - \,\mathbf{d}^\mathbf{p} ) + \frac{\mu}{2} \, || \mathbf{z}^\mathbf{p} - \sum_{\forall i\in \mathcal{N_M}} \mathbf{y}_i^\mathbf{p}   - \mathbf{d}^\mathbf{p}\,||^2_2
\end{split}
\label{aug_lagrangian}
\end{align}
where $\mu>0$ is a tuning parameter and $\boldsymbol{\lambda}$ is named \emph{dual variable} or \emph{shadow price}.
Being \eqref{aug_lagrangian}  not fully separable among the agents due to the presence of the quadratic term, the sequential iterative procedure described in \textbf{Algorithm 1} must be followed according to ADMM procedure.
 \begin{algorithm}
 \caption{Distributed economic dispatch through ADMM}
 \begin{algorithmic}
  \STATE {Define with $k$ the iteration index and initialize as: \linebreak $ k = 1 $ and $\boldsymbol{\lambda}^{0} = \boldsymbol{\lambda}_0$}
\vspace{1mm}
  \STATE {\textbf{while}   \textit{convergence} is not met} \vspace{1.5mm}
  \STATE {1). The MGCCs solve in parallel the following sub-problems using information about the previous iteration:}
\vspace{1.5mm}
\STATE{\hspace{4mm}\textbf{for all} $i \in \mathcal{N_M}$
\begin{align}
(\mathbf{x}_i^{\mathbf{p},k},\mathbf{y}_i^{\mathbf{p},k}) =  \underset{\text{s.t. \eqref{mg_constr}-\eqref{mg_output}}} {\text{argmin}} \, L_{\mu}(\mathbf{x}_{i}^{\mathbf{p}},\mathbf{y}_{i}^{\mathbf{p}},\boldsymbol{\lambda}^{k-1}, \mathbf{y}^{\mathbf{p},k-1}_{\forall j  \neq i},\mathbf{z}^{\mathbf{p},k-1})
\label{admm_mgs}
 \end{align}
\hspace{4mm}\textbf{end for} }
\vspace{1.5mm}
\STATE 2). The AGS gathers the optimal outputs of the MGs sub-problems, i.e. $\mathbf{y}_i^{\mathbf{p},k}$, and it solves the following sub-problem
\begin{align}
\mathbf{z}^{\mathbf{p},k} = \underset{\text{s.t. \eqref{AG_cmpct_AGconstr}}} {\text{argmin}} \; L_{\mu}(\mathbf{z}^{\mathbf{p}} ,\;\mathbf{y}^{\mathbf{p},k} _{\forall i},\;\boldsymbol{\lambda}^{k-1})
\label{admm_ags}
 \end{align}
\STATE 3). The dual variable $\boldsymbol{\lambda}$ is updated by the AGS based on the updated constraint residual
\begin{align}
\boldsymbol{\lambda}^{k} = \boldsymbol{\lambda}^{k-1} \; + \; \mu \cdot ( \mathbf{z}^{\mathbf{p},k}  - \sum_{\forall i\in \mathcal{N_M}} \mathbf{y}_i^{\mathbf{p},k}   -   \mathbf{d}^{\mathbf{p}} \,)
\label{admm_dual}
\end{align}
  \STATE {$k$ is updated for the next iteration:\quad $k=k+1$}
  \vspace{1.5mm}
    \STATE {\textbf{end while}}
 \end{algorithmic}
 \end{algorithm}
\\It could be pointed out that Step 1 of \textbf{Algorithm 1} does not allow to preserve the confidentiality since each $i$-th MGCC needs information about the optimal outputs of the other agents at the previous iteration, i.e. $ (\mathbf{y}^{\mathbf{p},k-1}_{\forall j  \neq i},\; \mathbf{z}^{\mathbf{p},k-1})$, in order to minimize $L_{\mu}$. However, defining the constraint residual as $\mathbf{r}^{\mathbf{p}} = \mathbf{z}^{\mathbf{p}}  - \sum_{\forall i\in \mathcal{N_M}} \mathbf{y}_i^{\mathbf{p}}  -   \mathbf{d}^{\mathbf{p}}$, it holds that
\begin{flalign}
\begin{split}
& L_{\mu}(\mathbf{x}_{i}^{\mathbf{p}},\mathbf{y}_{i}^{\mathbf{p}},\boldsymbol{\lambda}^{k-1}, \mathbf{y}^{\mathbf{p},k-1}_{\forall j  \neq i},\mathbf{z}^{\mathbf{p},k-1})=  \;f_i^p(\mathbf{x}_i^{\mathbf{p}}) \;-\; \boldsymbol{\lambda}^{k-1}{'} \mathbf{y}_i^{\mathbf{p}} \;+\\[3mm]&+ \;\frac{\mu}{2} ||\mathbf{z}^{\mathbf{p},k-1}   - \sum_{\substack{ \forall j \neq i \in \mathcal{N_M}}} \mathbf{y}^{\mathbf{p},k-1}_j  - \mathbf{y}_i^{\mathbf{p}} -  \mathbf{d}^{\mathbf{p}}||^2_2=\\&=f^p_i(\mathbf{x}_i^{\mathbf{p}}) - \boldsymbol{\lambda}^{k-1}{'} \mathbf{y}_i^{\mathbf{p}} - \frac{\mu}{2} || \;\mathbf{r}^{\mathbf{p},k-1} + \mathbf{y}_i^{\mathbf{p},k-1} - \mathbf{y}_i^{\mathbf{p}}\;||^2_2 \;=\\
&=L_{\mu}(\mathbf{x}_{i}^{\mathbf{p}},\mathbf{y}_{i}^{\mathbf{p}},\boldsymbol{\lambda}^{k-1},\mathbf{y}_{i}^{\mathbf{p},k-1}, \mathbf{r}^{\mathbf{p},k-1})
\label{admm_mgs_lagrangian}
\end{split}&&
\end{flalign}
Therefore, to perform Step 1 it is enough that the AGS provides each $i$-th MGCC with information about the dual variable and the coupling constraint residual. The MG$_i$ cost function to minimize is therefore defined by \eqref{admm_mgs_lagrangian}, and it includes the internal production costs $ f_i^p(\mathbf{x}_{i}^{\mathbf{p}})$ and some other terms  related to the output variables $\mathbf{y}_{i}^{\mathbf{p}}$. These terms were denoted in \eqref{muG_costfunction_compact} with the generic function $h_i^p(\mathbf{y}_{i}^{\mathbf{p}})$  and, at convergence, they express the cost/gain of the MG$_i$ for the provided  active power output and power reserves.
\\
The optimality of the asymptotic solution to \textbf{Algorithm 1} is guaranteed by the following proposition, which can be straightforwardly proven based on the proof in \cite[Section 3.4]{bertsekas1989parallel}.
\begin{prop}
The sequence $\{\mathbf{x}_{\forall i}^{\mathbf{p}},\mathbf{y}_{\forall i}^{\mathbf{p}},\mathbf{z}^{\mathbf{p}}\}^k$ generated by \textbf{Algorithm 1} is bounded and its limit points are in the set of the optimal solutions of the original problem \eqref{AG_cmpct}.		
\label{Prop2}
\end{prop}
Therefore, the sequence of the optimal value of the overall cost function, generated by \textbf{Algorithm 1} and defined as $\mathbf{J^{\textbf{\tiny{AG}},k}}=\sum_{\forall  i\in \mathcal{N_M}} f_i^p(\mathbf{x}_i^{\mathbf{p},k})  + g_z^p(\mathbf{z}^{\mathbf{p},k} ) $, converges to the optimal objective of the primal problem, i.e. $\mathbf{J^{\textbf{\tiny{AG}},k}} \longrightarrow\mathbf{J^{\textbf{\tiny{AG}},*}} $ as $k\rightarrow \infty$. 
Although this is an asymptotic result, the ADMM often converges in few tens of iterations with satisfactory accuracy, and the following termination criterion can be used:\linebreak $
 {||\mathbf{r}^{\mathbf{p},k}||}_2 \leq \epsilon_{r} \; {\wedge} \; ||\mathbf{z}^{\mathbf{p},k+1}-\mathbf{z}^{\mathbf{p},k}||_2 \leq \epsilon_{z}
$, where $\epsilon_{r}>0$ and $\epsilon_{z}>0$ are predefined tolerances \cite{boyd2011distributed}.\\	
Considering the practical implementation, the ADMM algorithm applied to the eAG management consists in the following iterative procedure: firstly the MGCCs perform in parallel their local optimization problems \eqref{admm_mgs}, based on the previous values of the constraint residual $\mathbf{r}^{\mathbf{p}}$, and of the internal shadow price $\boldsymbol{\lambda}$. Then, the AGS gathers the optimal values of the MGs output variables and solves its sub-problem \eqref{admm_ags}, considering the reserve minimum requirements and the external grid prices. Notice that the MGs internal information and optimization variables $\mathbf{x}_i^{\mathbf{p}}$ are not externally shared. Finally, the AGS updates the internal price $\boldsymbol{\lambda}$  based on the coupling constraint residual, as in \eqref{admm_dual}, which is communicated again to the MGs such that they can start again performing the first step of \textbf{Algorithm 1}. The dual variable $\boldsymbol{\lambda} = [\boldsymbol{\lambda^{r_p \uparrow}}\hspace{-1.8mm}{'},\boldsymbol{\lambda^{r_p \downarrow}}{'}, \boldsymbol{\lambda^p}{'}]{'}$ can be therefore interpreted as the vector of the internal negotiation prices between the AGS and the MGs, for the output active power ($\boldsymbol{\lambda^p}{'}$) and for the power reserves ($\boldsymbol{\lambda^{r_p \uparrow}}\hspace{-1.8mm}{'},\boldsymbol{\lambda^{r_p \downarrow}}{'})$. Differently from the external trade with the main utility where the selling and the buying prices are different, just one internal price exists for the MGs output active power. This is done on purpose since additional transmission-related fixed costs must not be considered internally to the eAG.
\section{Phase 2: power flow feasibility and reactive power planning}\label{sec:3}
A second optimization phase is necessary for twofold reason. First, to ensure that the active power flows computed by \textbf{Algorithm 1} are consistent with the constraints on nodal voltages and line currents. Second, to properly regulate the reactive power flows inside the eAG,  minimizing the power losses and ensuring the network electrical feasibility. This is made possible since MGs can have a significant role as reactive power producers/consumers, being clusters of several inverter-interfaced generation sources. 
As mentioned in the Introduction, since each MG is regarded as a single equivalent generator, this phase can be carried out by the AGS through a centralized approach, without compromising the computational feasibility and without violating the MGs privacy constraints.
The eAG network is assumed to be linked to the utility grid through a single point, referred as \textit{slack node} and conventionally numbered as the first node in $\mathcal{N}=\{1,\hdots,n\}$. The presence of an OLTC transformer on the slack node is also assumed, which adds a beneficial additional degree of freedom by allowing for the step-wise variation of the slack nodal voltage. This is often an effective solution to maintain the network voltages and currents in the proper range without affecting the scheduled active power set-points of the generation units \cite{farina2015model,casavola2017command}.
For clarity, from now on the optimal values of the variables computed in Section \ref{sec:2} by \textbf{Algorithm 1} will be denoted with the superscript *, e.g. $ \mathbf{p}_i^{\mathbf{\scriptstyle{mg,*}}}$.
\subsection{Modelling of the MGs as equivalent generators}\label{subsec:mg_reactive}
At this stage, the MG output reactive power coincides with the internal load demand since the generation units reactive power production/absorption has not been yet scheduled. The overall MG reactive power capability can be represented by aggregating the capabilities of each generation unit. For consistency with Section \ref{sec:2}, the overall MG reactive power capability is expressed through up/down reserves. The following relations hold:
\begin{align}
\mathbf{q}_{\,i}^\mathbf{\scriptstyle{mg,*}} \;&=\; \mathbf{q}_{\,i}^\mathbf{\scriptstyle{l}} \label{mg_q_I}
\\
 \mathbf{r}_{\mathbf{q},i}^\mathbf{\scriptstyle{mg}\uparrow\scriptstyle{*}} &=  \sum_{ j_i=1}^{n^g_i} \mathbf{\bar{q}}_{j_i}^\mathbf{g}(\mathbf{p}_{j_i}^\mathbf{g,*}) +\sum_{k_i=1}^{n^b_i} \mathbf{\bar{q}^b}_{k_i}+ \sum_{p_i=1}^{n^r_i}\mathbf{\bar{q}^r}_{p_i}
\\
\mathbf{r}_{\mathbf{q},i}^\mathbf{\scriptstyle{mg}\downarrow\scriptstyle{*}}  &= - \sum_{ j_i=1}^{n^g_i} \mathbf{\underline{q}}_{j_i}^\mathbf{g}(\mathbf{p}_{j_i}^\mathbf{g,*}) -\sum_{k_i=1}^{n^b_i} \mathbf{\underline{q}}_{k_i}^\mathbf{b}- \sum_{p_i=1}^{n^r_i}\mathbf{\underline{q}}_{p_i}^\mathbf{r}\label{mg_q_rdwn_I}
\end{align}
where the generators reactive power limits are expressed as functions of the active power production, as usual for fuel-based generators, \cite{adibi1994reactive}. The MG reactive power output and reserves defined in \eqref{mg_q_I}-\eqref{mg_q_rdwn_I}  have been also denoted with the superscript $*$ since they are not optimization variables in this phase, but their values are computed based on the outcomes of \emph{Phase 1}. Actually, prior to perform  \emph{Phase 2}, MGs must communicate to the AGS not only their optimal active power trends, but also the active and reactive reserve capabilities.
The reactive power variation $\mathbf{\Delta q}^{\mathbf{\scriptstyle{mg}}}$ with respect to $\mathbf{q}^{\mathbf{\scriptstyle{mg,*}}}$ is then scheduled in \emph{Phase 2}, and this must be bounded by the up/down reserves:
\begin{align}
-\mathbf{r}_{\mathbf{q},\,i}^{\mathbf{\scriptstyle{mg}\downarrow\scriptstyle{*}}} \; \leq \;  \mathbf{\Delta q}_{i}^{\mathbf{\scriptstyle{mg}}}\; \leq \; \mathbf{r}_{\mathbf{q},\,i}^{\mathbf{\scriptstyle{mg}\uparrow\scriptstyle{*}}}
\label{reactive_power_var_lim}
\end{align}
It would be desirable that the optimal active power trends achieved by \textbf{Algorithm 1} are not varied, since it would be sub-optimal with respect to problem \eqref{AG_cmpct}. To this regard, since the MG active power output has been locally constrained, see \eqref{mg_output_constr}, it is expected that just small active power variations may be needed, e.g. in case more MGs are injecting an excessive amount of active power in the same line leading to over-current issues. Also in this case, the active power reserves scheduled in \emph{Phase 1} can be used to limit the active power variations of each MG. Therefore, denoting with $\Delta \mathbf{p}_{i}^{\mathbf{\scriptstyle{mg}}}$ the active power variation with respect to $\mathbf{p}_{\,i}^\mathbf{\scriptstyle{mg,*}}$, it holds that:
\begin{align}
-\mathbf{r}_{\mathbf{p},\,i}^{\mathbf{\scriptstyle{mg}\downarrow\scriptstyle{,*}}} \; \leq \;  \mathbf{\Delta p}_{i}^{\mathbf{\scriptstyle{mg}}}\; \leq \; \mathbf{r}_{\mathbf{p},\,i}^{\mathbf{\scriptstyle{mg}\uparrow\scriptstyle{,*}}}
\label{active_power_var_constr}
\end{align}
Nevertheless, the variation of the MGs active power outputs should not compromise the minimum eAG frequency reserve requirements \eqref{AG_compact_centr_reserve_min} defined in Section \ref{sec:2}. Therefore the following additional constraints must be considered:
\begin{align}
& \sum_{\forall i \in \mathcal{N_M}} ( \mathbf{r}_{\mathbf{p},\,i}^{\mathbf{\scriptstyle{mg}\uparrow}\scriptstyle{,*}}-\Delta \mathbf{p}_{i}^{\mathbf{\scriptstyle{mg}}}) \; \geq \; \mathbf{\underline{r}_p^{\textbf{\tiny{AG}}\uparrow}}  \label{up_reserve_reactive}\\
&\sum_{\forall i \in \mathcal{N_M}} (\mathbf{r}_{\mathbf{p},\,i}^{\mathbf{\scriptstyle{mg}\downarrow}\scriptstyle{,*}}+\Delta \mathbf{p}_{i}^{\mathbf{\scriptstyle{mg}}}) +\sum_{\forall j_r \in \mathcal{N_R}} \mathbf{p}^{\textbf{\tiny{R}}}_{j_r} \; \geq \; \mathbf{\underline{r}_p^{\textbf{\tiny{AG}}\downarrow}}\label{down_reserve_reactive}
\end{align}
To sum up, the MGs are modelled as PQ generation nodes with predefined active and reactive power trends, i.e. $ \mathbf{p}^{\mathbf{\scriptstyle{mg,*}}}$ and $ \mathbf{q}^{\mathbf{\scriptstyle{mg,*}}}$, allowing for active and reactive power set-points variations, $ \Delta \mathbf{p}^{\mathbf{\scriptstyle{mg,*}}}$ and $\Delta \mathbf{q}^{\mathbf{\scriptstyle{mg,*}}}$ which must respect predefined capability limits, defined by   \eqref{reactive_power_var_lim}-\eqref{down_reserve_reactive}.\\
It should be considered that some approximations are introduced by the MG equivalent representation. For instance, if the active power variation of certain MG  implies that a battery power production must vary at a precise time instant, the corresponding reserve for the next time steps could change as well since the SOC  changes, see \eqref{batt_upreserve_costr_2}-\eqref{batt_downreserve_costr_2}. However, as already mentioned, consistent active power variations are not expected and these will be highly minimized in the optimization problem of \emph{Phase 2} acting as slack variables, used just if electrical feasibility is compromised.
%
\subsection{OLTC model}
The on-load tap-changer can assume a finite number of positions allowing to vary the slack voltage around its nominal value. The slack voltage $V_1(t)$ is therefore modelled to take $2\, n_s$ discrete values, i.e. $V_1(t)=V^o\pm k\Delta {V}^s$ where $k=0,\dots,n_s$ and $n_s=\lfloor  {\Delta \bar{V}^s}/{\Delta {V}^s}  \rfloor$, see Table~\ref{tab:MG_active_opt_var}.\\ In view of this, the OLTC daily voltage profile $\mathbf{V_1}$ is formally defined as
%
\begin{equation}
\mathbf{V}_1 = V^o\,( \,\mathbf{1}_{\scriptstyle{N}} \; + \; \sum_{k=-n_s}^{{n}_s} \boldsymbol{\theta}_{k}^\mathbf{s}\, k
\Delta V^s) \label{slack_OLTC}\end{equation}
where $\mathbf{1}_N$ is the $N$-dimensional unity vector and $\boldsymbol{\theta}_{k}^\mathbf{s}$ (for $k=-n_s,\dots,n_s$) are vectors of optimization variables whose entries take values in $\{0,1\}$, defining when a certain position $k$ is selecter over the day. Moreover, the following constrain is stated
\begin{equation}
\sum_{k=-n_s}^{{n}_s}\boldsymbol{\theta}_{k}^\mathbf{s} \;=\;\mathbf{1}_{\scriptstyle{N}}
\label{sum_steps_1}
\end{equation}
in such a way that, at each time step of the day, only one OLTC position can be selected.
\subsection{Network model}
The power flow equations must be considered to model the network electrical variables. They are nonlinear functions defining the nodal powers based on the network voltages and on the \textit{nodal admittance matrix}. Moreover, combining these equations, also some other variables can be calculated such as the active power losses and the magnitude of the flowing current for each line of the network. These equations are not explicitly expressed here as they be easily recovered from the literature, see \cite{marconato2002electric}. The following generic notation is therefore used
\begin{align}
\mathbf{P}_j\; &= \; \boldsymbol{f}_j^{\boldsymbol{P}} (\mathbf{V}_{1,\hdots,n},\boldsymbol{\delta}_{1,\hdots,n},Y)	\quad && \forall j \in \mathcal{N}	\label{P_pf}\\
\mathbf{Q}_j \; &= \; \boldsymbol{f}_j^{\boldsymbol{Q}} (\mathbf{V}_{1,\hdots,n},\boldsymbol{\delta}_{1,\hdots,n},Y)		\quad &&\forall j \in \mathcal{N}		\label{Q_pf}\\
\mathbf{P}_{\scriptstyle{i,j}}^\mathbf{loss} \; &= \; \boldsymbol{f}_{\scriptstyle{i,j}}^{\boldsymbol{loss}} (\mathbf{V}_i,\boldsymbol{\delta}_i,\mathbf{V}_j,\boldsymbol{\delta}_j,Y) \quad &&\forall (i,j) \in \mathcal{E}\\
\mathbf{I}_{\scriptstyle{i,j}} \; &= \; \boldsymbol{f}_{\scriptstyle{i,j}}^{\boldsymbol{I}} (\mathbf{V}_i,\boldsymbol{\delta}_i,\mathbf{V}_j,\boldsymbol{\delta}_j,Y) \quad &&\forall (i,j) \in \mathcal{E}\label{I_pf}
\end{align}
where $\boldsymbol{f}_j^{\boldsymbol{P}} $, $\boldsymbol{f}_j^{\boldsymbol{Q}} $, $\boldsymbol{f}_{\scriptstyle{i,j}}^{\boldsymbol{loss}}$ and $\boldsymbol{f}_{\scriptstyle{i,j}}^{\boldsymbol{I}} $ are vectors of static nonlinear functions expressing the $j$-th nodal active power, the $j$-th nodal reactive power, the $(i,j)$ line active power loss and the $(i,j)$ line current magnitude over the whole optimization horizon, respectively.
For each node of the eAG network, also the power balance for each node must be stated, linking the nodal powers defined in \eqref{P_pf}-\eqref{Q_pf} with the output powers of the eAG elements; it follows that $\forall j \in \mathcal{N}$
\begin{align}
\begin{split}
\mathbf{P}_j = \hspace{-2mm}&\sum_{\forall i \in \mathcal{N_M}} \sigma^M_{i,j} \,(\mathbf{p}_{i}^{\mathbf{\scriptstyle{mg,*}}} + \mathbf{\Delta p}_{i}^{\mathbf{\scriptstyle{mg}}})+\hspace{-2mm} \sum_{\forall j_l \in \mathcal{N_L}}\sigma^L_{j_l,j}\, \mathbf{p}_{j_l}^{\mathbf{L}}\;+\\+ &\sum_{\forall j_r \in \mathcal{N_R}}\sigma^R_{j_r,j}\, \mathbf{p}_{j_r}^{\mathbf{R}}\label{nodal_bal_p}
\end{split}\\[2mm]
\begin{split}
\mathbf{Q}_j =\hspace{-2mm} &\sum_{\forall i \in \mathcal{N_M}} \sigma^M_{i,j} \,(\mathbf{q}_{i}^{\mathbf{\scriptstyle{mg,*}}} + \mathbf{\Delta q}_{i}^{\mathbf{\scriptstyle{mg}}})+\hspace{-2mm} \sum_{\forall j_l \in \mathcal{N_L}}\sigma^L_{j_l,j}\, \mathbf{q}_{j_l}^{\mathbf{L}}\label{nodal_bal_q}
\end{split}
\end{align}
where  $\sigma^{M,R,L}_{\alpha,\beta}$ are boolean scalars defined to be equal to 1 just in case the $\alpha$-th element is connected to the $\beta$-th node.
Finally, the electrical constraints must be also introduced such that nodal voltages and line currents are kept inside the allowed range
\begin{align}
&\underline{\mathbf{V}}_j \; \leq \; \mathbf{V}_j \; \leq \; \bar{\mathbf{V}}_j \qquad  &&\forall \, j \in \mathcal{N} \label{V_constr}\\
&\mathbf{I}_{i,j}(t) \; \leq \; \bar{\mathbf{I}}_{i,j} \qquad  \;\,&& \forall \,(i,j) \in \mathcal{E}&&\label{I_constr}
\end{align}
\subsection{Phase 2 optimization problem}
The overall optimization problem to be solved in \emph{Phase 2} can now be stated:
\begin{subequations}
\begin{flalign}
\min_{\substack{\Delta \mathbf{p}_{\scriptstyle{\forall i }}^{\mathbf{\scriptstyle{mg}}}, \Delta \mathbf{q}_{\scriptstyle{\forall i }}^{\mathbf{\scriptstyle{mg}}}\\\boldsymbol{\theta}_{\forall i_s}^{\mathbf{s}}}} \Big \{&\sum_{\forall(i,j)\in \mathcal{E}} || \mathbf{P}_{i,j}^\mathbf{\scriptstyle{loss}}||^{\scriptstyle{2}}_{\scriptstyle{2}}+\hspace{-2mm} \sum_{\forall  i \in \mathcal{N_M}} \hspace{-1mm}\zeta_i\, ||{\Delta \mathbf{p}_{\scriptstyle{i}}^{\mathbf{\scriptstyle{mg}}}}||^{\scriptstyle{2}}_{\scriptstyle{2}}\,+  \nonumber &&\\[-3mm]&\hspace{2mm}+\sum_{k=-n_s}^{{n}_s}   w_s k^2	\,||\boldsymbol{\theta}_{k}^{\mathbf{s}}||^{\scriptstyle{2}}_{\scriptstyle{2}}\; \Big \}&&
\label{opf_cost_function}
\end{flalign}
$subject \;to$
\begin{equation}
\left. \hspace{-5mm}
\begin{aligned}
&\textit{MGs active and reactive power limits}:
  &\text{\eqref{reactive_power_var_lim},\,\eqref{active_power_var_constr} } \hspace{-1mm}\\
&\textit{Active reserve requirements}:  &\text{\eqref{up_reserve_reactive},\,\eqref{down_reserve_reactive}}  \\
&\textit{OLTC slack voltage modelling}:  &\text{\eqref{slack_OLTC},\,\eqref{sum_steps_1}}  \\
& \textit{Power flow equations}:  &\text{\eqref{P_pf}-\eqref{I_pf}}  \\
& \textit{Nodal power balances}:  &\text{\eqref{nodal_bal_p},\,\eqref{nodal_bal_q}}  \\
& \textit{Voltages and currents limits}:  & \text{\eqref{V_constr},\,\eqref{I_constr}}
\end{aligned}
\hspace{1mm} \right \} \hspace{-3mm}
\end{equation}
\label{centr_opf}
\end{subequations}
\\
%
%
%
The first term in \eqref{opf_cost_function} penalizes the active power losses. The second one minimizes the MG active power variations through the  parameter $\zeta_i>0$, set to a very high value such that the active power variations are allowed just if the problem feasibility is compromised. The third term penalizes the variations of the OLTC position. The corresponding weight, $w_sk^2>0$, is defined such that the larger is the position variation, the higher is the cost.
The stated optimization problem is mixed-integer and non-convex, and it can be centrally solved by means of global optimization techniques \cite{burer2012non}.
\section{Phase 3: Final scheduling of the MGs generation units} \label{sec:4}
A final optimization procedure is needed to schedule the power references of the MGs generation units consistently with the active and reactive optimal profiles obtained in Phases 1 and 2. 
Note that, in the computation of the latter,  some approximations were introduced by modelling the MGs as equivalent generators. Therefore,  Phase 3 is necessary also to verify that the MGs aggregation still provides the minimum required power reserves, even though the constraints \eqref{up_reserve_reactive}-\eqref{down_reserve_reactive} have been  considered in Phase 2. Denoting the optimal values of the variables computed in Phase 2 with the superscript $**$, the following optimization problem must be solved.
\begin{subequations}
\begin{align}
&\min_{\substack{ \mathbf{x}_{\forall i}^{\mathbf{p}},\mathbf{y}_{\forall i}^{\mathbf{p}} \\ \mathbf{x}_{\forall i}^{\mathbf{q}},\mathbf{y}_{\forall i}^{\mathbf{q}}\\ \mathbf{z}^\mathbf{p} }} \Big  \{ \sum_{\forall  i\in \mathcal{N_M}} \hspace{-2mm}\Big [ {f}_i^{{p}} (\mathbf{x}_i^{\mathbf{p}}) +{f}_i^{{q}} (\mathbf{x}_i^{\mathbf{q}})+\gamma_i^p\, ||\boldsymbol{\epsilon}_i^{\boldsymbol{p}}||_{\scriptstyle{2}}^{\scriptstyle{2}} +\gamma_i^q \,||\boldsymbol{\epsilon}_i^{\boldsymbol{q}}||_{\scriptstyle{2}}^{\scriptstyle{2}}    \Big ]\Big \}
\label{AG_cmpct_globalcf_2}
\end{align}
\\[-5mm]subject to\\
\begin{empheq}[right=\empheqrbrace]{align}
\begin{split}
& \hspace{-15mm} A_i^p\,
\mathbf{x}_i^\mathbf{p}\, \leq\, b_i^p \\
 &\hspace{-15mm}\,\,\mathbf{y}_i^\mathbf{p}\, = \, C_i \,\mathbf{x}_i^\mathbf{p} \, + \, M_i\, \mathbf{d}_{i}^\mathbf{p}\label{mg_active_constr_2}
\end{split}
 \\[2mm]
 \begin{split}
 &\hspace{-15mm}A_i^q\,
\mathbf{x}_i^\mathbf{q}\, \leq\, b_i^q(\mathbf{x}_i^\mathbf{p}) \\
 & \hspace{-15mm} \,\,\mathbf{y}_i^\mathbf{q}\, = \, C_i^q \,\mathbf{x}_i^\mathbf{q} \, + \, M_i^q\, \mathbf{d}_{i}^\mathbf{q}\label{mg_reactive_constr_2}
\end{split}
\\[2mm]
 \begin{split}
 &\hspace{-15mm}\; \mathbf{p}_{\,i}^\mathbf{\scriptstyle{mg}} \; = \; \mathbf{p}_{i}^\mathbf{\scriptstyle{mg,*}} + \mathbf{\Delta p}_{i}^\mathbf{\scriptstyle{mg,**}} \; + \boldsymbol{\epsilon}_i^{\boldsymbol{p}} \; \; \\
 & \hspace{-15mm}\; \mathbf{q}_{i}^\mathbf{\scriptstyle{mg}} \; = \; \mathbf{q}_{i}^\mathbf{\scriptstyle{mg,*}} + \mathbf{\Delta q}_{i}^\mathbf{\scriptstyle{mg,**}} \; + \boldsymbol{\epsilon}_i^{\boldsymbol{q}}\label{mg_stp_2}
\end{split}
\end{empheq}\\[-5mm]   $ \hspace*{64mm}  \forall i \in \mathcal{N_M} $\\[-6mm]
 \begin{align}
 &\hspace{-25mm} \mathbf{z}^{\mathbf{p}}   \; \in \mathcal{Z}^p \label{AG_constr_2}
 \\
 &\hspace{-25mm} \mathbf{z}^{\mathbf{p}} \;- \sum_{\forall i\in \mathcal{N_M}} \mathbf{y}_i^{\mathbf{p}} \;  - \,\mathbf{d}^{\mathbf{p}}  \;= \;\mathbf{0}\label{global_active_constr_2}
 \end{align} \label{rescheduling_opt}
\end{subequations}
Following the same reasoning of Section 		\ref{sec:2}, the internal and output reactive power variables are compactly defined as $\mathbf{x}_i^{\mathbf{q}}$ and $\mathbf{y}_i^{\mathbf{q}}$, while the MGs internal reactive power constraints are given by \eqref{mg_reactive_constr_2}.  The constraints \eqref{mg_stp_2} are here introduced to make the MGs follow the active and reactive power outputs dictated by the optimization problems in Phases 1 and 2. Two slack variables are used in \eqref{mg_stp_2}, i.e. $\boldsymbol{\epsilon}_i^{\boldsymbol{p}}$ and $\boldsymbol{\epsilon}_i^{\boldsymbol{q}}$, to avoid infeasibility problems, and they are  minimized in the cost function \eqref{AG_cmpct_globalcf_2} through  the weights $\gamma_i^p>0$ and $\gamma_i^q>0$, taking very large values. 
Considering the cost function, the MG productions cost are included through the function ${f}_i^{{p}}$, already defined in Section \ref{sec:2}, while a function ${f}_i^{{q}}$ is also introduced to properly distribute the reactive power burden among the MG generation units.
The active reserve minimum requirements, defined in \eqref{AG_constr_2} as explained in Section \ref{sec:2}, and the coupling constraint \eqref{global_active_constr_2}, are also included in the problem formulation of Phase 3.
\\
The optimization problem \eqref{rescheduling_opt} cannot be solved in a centralized way because this would require the MGs to externally share their internal information. Therefore, two different approaches must be used based on the outcome of Phase 2: 
\begin{itemize}
\item \emph{Phase 3(a)}: No active power variations have been requested in Phase 2, meaning that the power trends scheduled in Phase 1 did not lead to over-current or over-voltage issues. This means that the reserve provision requirements, i.e. \eqref{AG_constr_2}, is already respected using the active power references scheduled in Phase 1. Therefore,  \eqref{rescheduling_opt} can be solved in a fully decentralized fashion neglecting constraints \eqref{AG_constr_2} and \eqref{global_active_constr_2}, i.e. each MG solves its local sub-problem scheduling the reactive power references of its generation units, keeping as active power set-points the ones computed in Phase 1.
\item \emph{Phase 3(b)}: In case active power variations have been requested, a distributed algorithm must be used in order to ensure that constraint \eqref{AG_constr_2} is globally respected. The same procedure of \textbf{Algorithm 1} can be used in this phase, by adopting the following Augumented Lagrangian function\\
\begin{align*}
\begin{split}
& \hspace{-5mm} \tilde{L}_{\mu}(\mathbf{x}_{\forall i}^\mathbf{p},\mathbf{y}_{\forall i}^\mathbf{p},\mathbf{x}_{\forall i}^\mathbf{q},\mathbf{y}_{\forall i}^\mathbf{q}, \mathbf{z}^\mathbf{p} ,\boldsymbol{\tilde{\lambda}})\, = \sum_{\forall  i\in \mathcal{N_M}} \hspace{-2mm}\Big [ \, \boldsymbol{f}_i^{\boldsymbol{p}} (\mathbf{x}_i^{\mathbf{p}}) \,+\,\boldsymbol{f}_i^{\boldsymbol{q}} (\mathbf{x}_i^{\mathbf{q}}) \,+\\& \hspace{-8mm}+ \;\gamma_i^p\, ||\boldsymbol{\epsilon}_i^{\boldsymbol{p}}||_{\scriptstyle{2}}^{\scriptstyle{2}} \,+\,\gamma_i^q \,||\boldsymbol{\epsilon}_i^{\boldsymbol{q}}||_{\scriptstyle{2}}^{\scriptstyle{2}}   \; \Big ] \;+\; \boldsymbol{\tilde{\lambda}}'(\mathbf{z}^\mathbf{p} - \sum_{\forall i\in \mathcal{N_M}} \mathbf{y}_i^\mathbf{p}   - \mathbf{d}^\mathbf{p}) \;+\\&\hspace{-8mm}+  \frac{\mu}{2} \, || \mathbf{z}^\mathbf{p} - \sum_{\forall i\in \mathcal{N_M}} \mathbf{y}_i^\mathbf{p}   - \mathbf{d}^\mathbf{p}\,||^2_2 
\end{split}
\end{align*}
\\
where $\boldsymbol{\tilde{\lambda}}$ is the dual variable. In this way, the active and reactive power references of each generation source are computed consistently with the requirements of Phase 1 and Phase 2. 
\end{itemize}
\subsection*{Final considerations}
It could be argued that separating the whole scheduling problem in three different phases may be not necessary and it could lead to a sub-optimal solution. Actually, the problem constraints and cost functions of Phase 1, 2 and 3 can be clustered in a unique non-convex optimization problem which, in principle, could be solved through a distributed algorithm.\\
However, it has been chosen to adopt this three steps procedure for two main reasons:

\begin{enumerate}
\item Applying distributed optimization algorithms to non-convex problems leads to sub-optimal solution, and not-generalized convergence results \cite{boyd2011distributed}. Possible solution algorithms exist but they usually have a complex structure, e.g. in case agents have to communicate the gradients and Hessians of their local cost functions and constraints \cite{houska2016augmented}, or they have a quite slow convergence \cite{magnusson2015distributed}.
\item The proposed structure allows to perform the optimization phases in different time instants and at different time rates. For instance, Phase 2 and 3 can be periodically performed in the online management of the MGs, using updated information about the load and renewable power trends, ensuring the MGs support for the voltage and current regulation, and for the minimization of power losses. The same holds for Phase 1, which may be used to reschedule the eAG active power production in case the energy prices or the reserve requirements are varied over the day.
\end{enumerate}
\section{Numerical results}\label{sec:5}
The simulations have been carried out considering two case studies: the IEEE 13-bus and the IEEE 37-bus system networks with a varying number of MGs; the network characteristics are derived from \cite{kersting2001radial}.  Morever, for the IEEE 13-bus system, the presence of an OLTC has been supposed at the interface with the main grid, while for the IEEE 37-bus system a direct connection with the main grid is modelled.
\subsection{Case study 1: IEEE 13-bus system}\label{sec:51}
\begin{figure}[t!]
\centering
\includegraphics[width=0.9\linewidth]{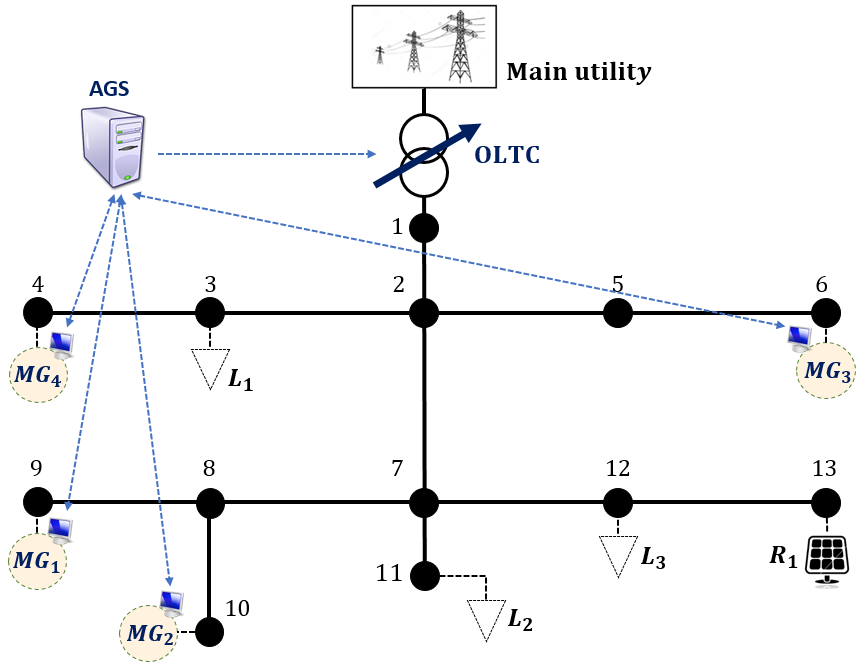}
\caption{Case study 1: eAG network topology: IEEE 13 bus system}
\label{fig:ieee13}
\end{figure}

\begin{figure}[!t]
\centering
\subfloat[]{\includegraphics[width=0.5\linewidth]{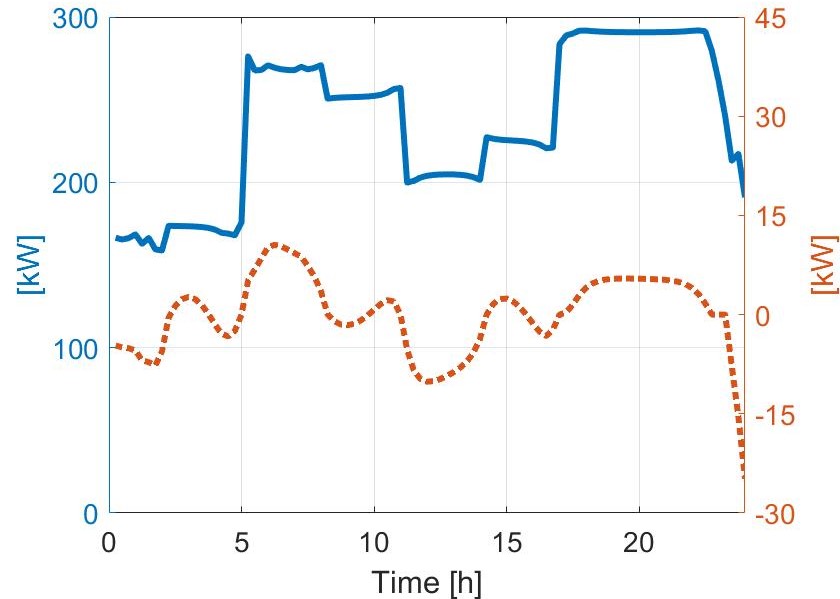}}
\subfloat[]{\includegraphics[width=0.5\linewidth]{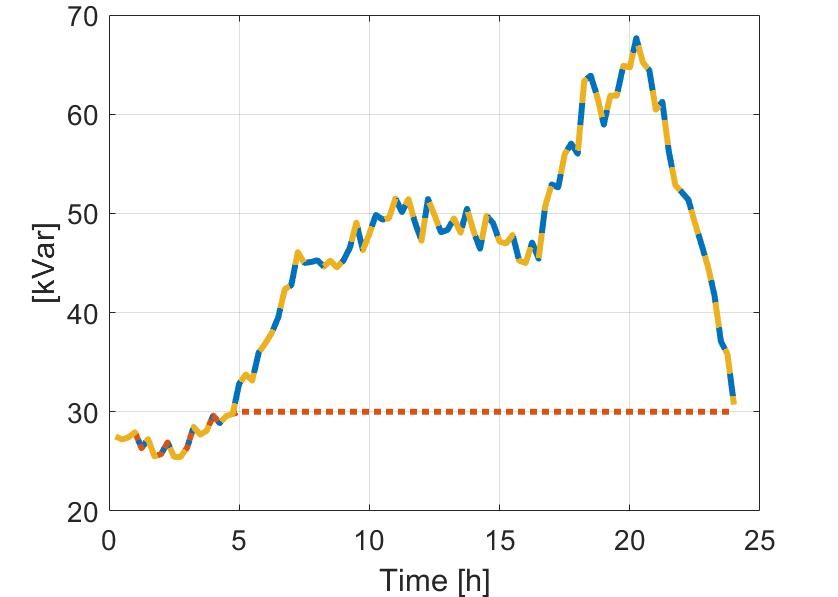}}
\caption{Case study 1: (a) MG2 active power set-points: Generator (solid), battery  (dotted); (b)  MG2 reactive power set-points: Micro-generator (solid), battery  (dotted), renewable source (dashed).}
\label{fig:REsch}
\end{figure}

 As shown in Figure~\ref{fig:ieee13}, an aggregation of four MGs (MG1, $\hdots$ , MG4), three loads (L1, L2, L3) and a non-controllable generation source (R1) is initially considered. The OLTC transformer is characterized by a maximum voltage deviation $\Delta \bar{V}=10 \%$ and by voltage step-size $\Delta {V}=2 \%$.
 In Table \ref{tab:MGs_gen}, the MGs generation units characteristics are described, where the abbreviation "Gen." stands for fuel-based generator, "Batt." for battery and "Res." for renewable energy source. For the following numerical results, the storage systems are all supposed to start with $SOC = 50\%$ and are limited between $20\%$ and $80\%$. The used active power trends for the load and renewable sources are shown in Figure \ref{fig:forecasts}(a)-(c)-(d). It is assumed that the loads reactive power trends are all characterized by a constant power factor of $0.8$. In Figure \ref{fig:forecasts}(b), the day-ahead energy selling and buying price trends are shown, while, for the sake of simplicity, the up and down power reserve prices are kept equal to a constant value for the whole day, i.e.  $\rho^{\uparrow}_{r_p}=\rho^{\downarrow}_{r_p}= 0.004 $ \euro /kWh. Concerning the ancillary services provision, the whole eAG must provide at least 100 kW of up and down power reserve at each instant of the day, while the line currents are limited by $170\, A$ and the nodal voltages by a maximum variation of $10\%$.
\begin{table}
\caption {Case study 1: MGs generation units}
\label{tab:MGs_gen}
\centering
\begin{tabular}{c c c c c}
\hline
&  & & &\\[-0.7em]
 Owner (Type) & $(\underline{p},\,\bar{p})$ & $(\underline{q},\,\bar{q})$ & $C^b$ & Costs  \\ \hline
MG1 (Gen.) & {(20, 250)} & $\pm\, \sqrt{(280 \text{\,-\,} p^g)}$ & $-$ & (0.02, 6, 30)$e^{-3}$ \\           MG1 (Gen.) & (20, 250) & $ \pm\, \sqrt{(280 \text{\,-\,} p^g)}$ & $-$ & (0.03, 6, 30)$e^{-3}$ \\
\vspace{1mm}MG1 (Res.) & (0, 80) & $\pm $80  & $-$ & $-$ \\
MG2 (Gen.) & (20, 300) & $\pm\, \sqrt{(330 \text{\,-\,} p^g)}$ & $-$ & (0.04, 1, 30)$e^{-3}$ \\
MG2 (Batt.) & $\pm$ 30 & $\pm$ 30 & 40 & 0.95$e^{-4}$  \\
\vspace{1mm}MG2 (Res.) & (0, 150) & $\pm 100 $ & $-$ & $-$ \\
MG3 (Gen.) & (10, 80) & $ \pm\, \sqrt{(90 \text{\,-\,} p^g)}$ & $-$ & (0.001, 3, 30)$e^{-3}$ \\
\vspace{1mm}MG3 (Res.) & (0, 30) & $\pm$ 30  & $-$ & $-$ \\
MG4 (Batt.) & $\pm$ 40 & $\pm$ 40 & $50$ & 1.15$e^{-4}$ \\
MG4 (Batt.) & $\pm$ 30 & $\pm$ 30 & $40$ & 1.10$e^{-4}$ \\
MG4 (Res.) & (0, 60) & $\pm$ 50  & $-$ & $-$ \\
 \hline
\end{tabular}
\end{table}
\begin{figure}[!t]
\centering
\subfloat[]{\includegraphics[width=0.5\linewidth]{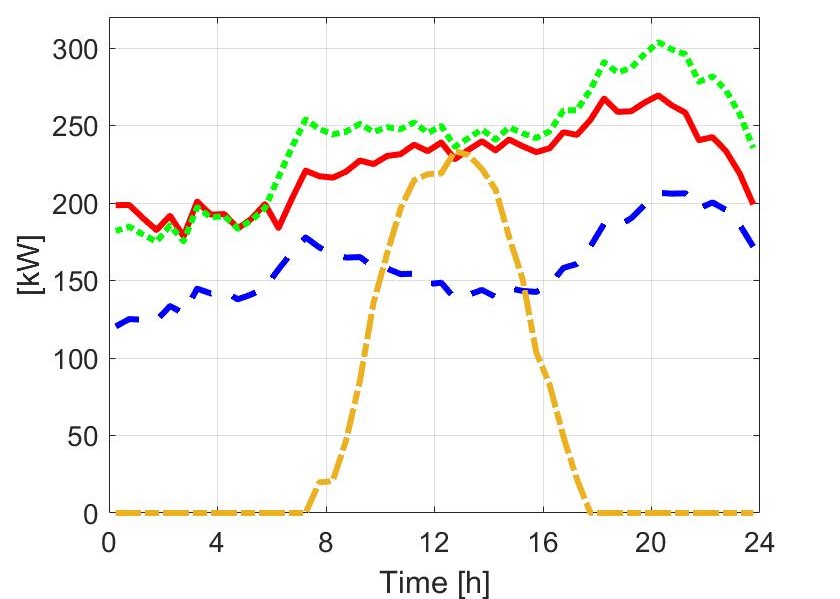}}
\subfloat[]{\includegraphics[width=0.5\linewidth]{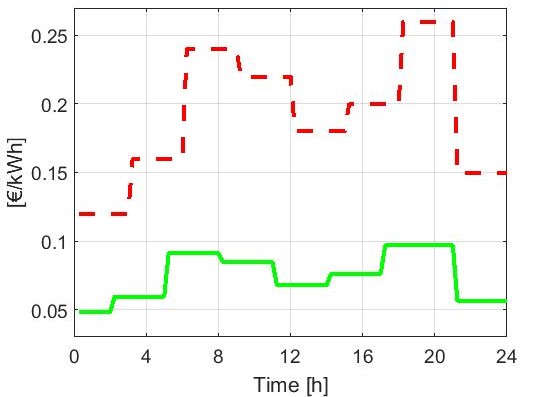}} \\[-3mm]
\subfloat[]{\includegraphics[width=0.5\linewidth]{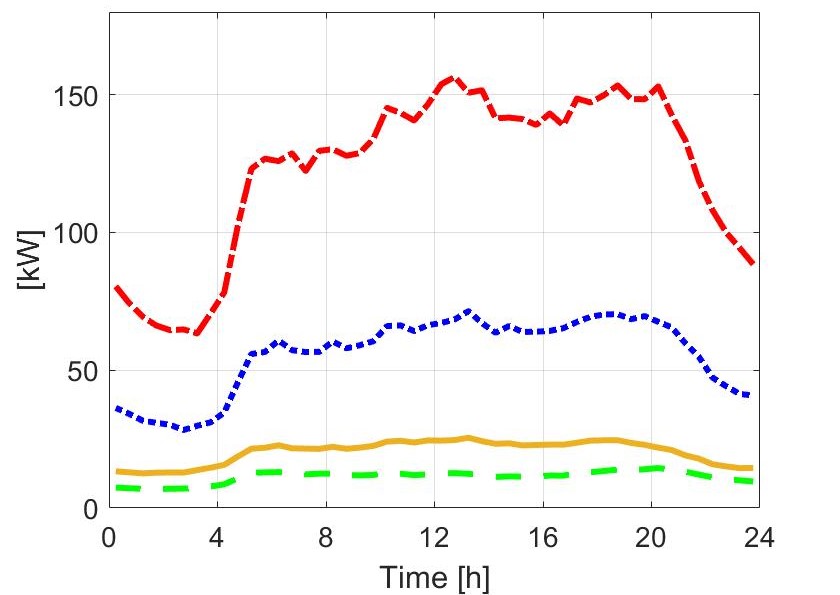}}
\subfloat[]{\includegraphics[width=0.5\linewidth]{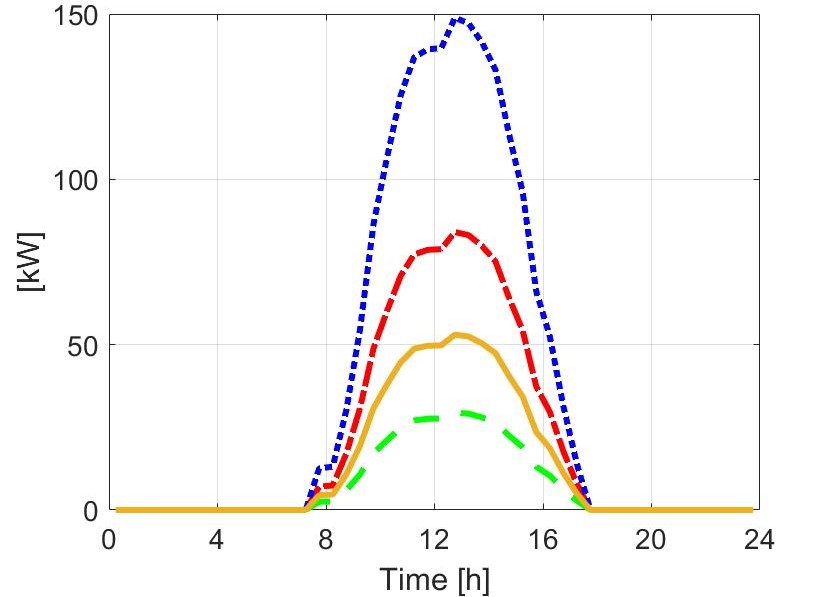}}
\caption{Case study 1: (a) eAG active power forecasts: L1(solid) L2(dashed) and L3(dotted) absorption, and R1(dashed-dotted) production; (b) Day-ahead energy prices: buying price (dashed), selling price (solid); (c)-(d)  MG1(dashed-dotted), MG2(dotted), MG3(dashed) and MG4(solid) load active power and renewable source production.}
\label{fig:forecasts}
\end{figure}
\\
\textbf{Algorithm 1} has been implemented in MATLAB and the corresponding optimization problems are solved using IBM ILOG CPLEX Optimization Studio. For the considered aggregation case, the ADMM algorithm converges to the global optimum in about 100 iterations.
Figures \ref{fig:ApoutRup}(a)-(b) show the scheduled eAG output power and provided up-reserve trends, where it is evident that the minimum required amount of power reserve is always globally provided. Considering the dual variables defined in \eqref{aug_lagrangian}, their steady-state solution can be interpreted as negotiation prices between the AGS and the MGs. Indeed, the internal energy price, defined by $\boldsymbol{\lambda^p}/\tau$, converges to the selling price when the eAG exports power, and to the buying price when the eAG absorbs power; in case the eAG does not exchange energy with the main utility, the internal energy price takes an intermediate optimal value, see Figure \ref{fig:ApoutRup}(c).  On the other hand, the internal price for the provided up reserve, $\boldsymbol{\lambda^{r_p \uparrow}}\hspace{-1.8mm}/\tau$, generally converges to the main utility reserve price, apart from the last part of the day where it increases to force the MGs to respect the minimum required amount, as depicted in Figure \ref{fig:ApoutRup}(d).
\begin{figure}[!t]
\centering
\subfloat[]{\includegraphics[width=0.5\linewidth]{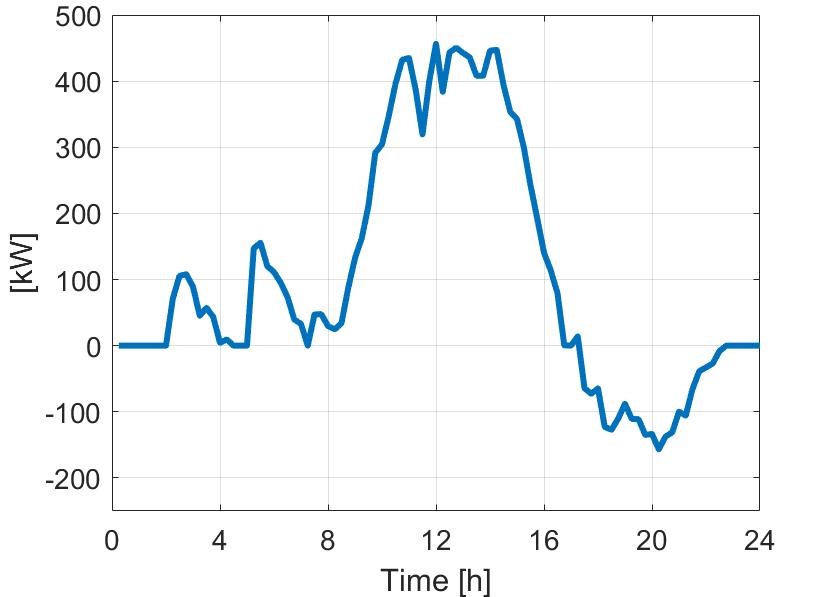}}
\subfloat[]{\includegraphics[width=0.5\linewidth]{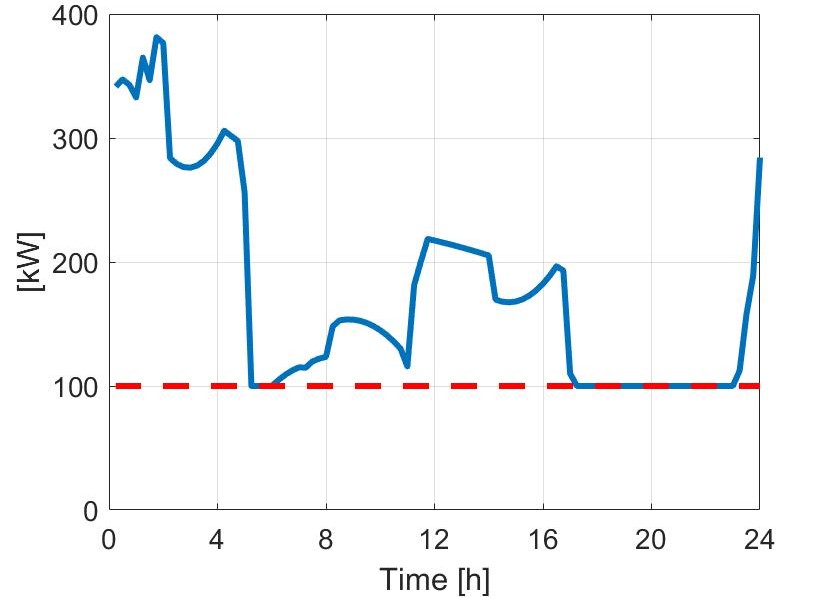}} \\[-3mm]
\subfloat[]{\includegraphics[width=0.5\linewidth]{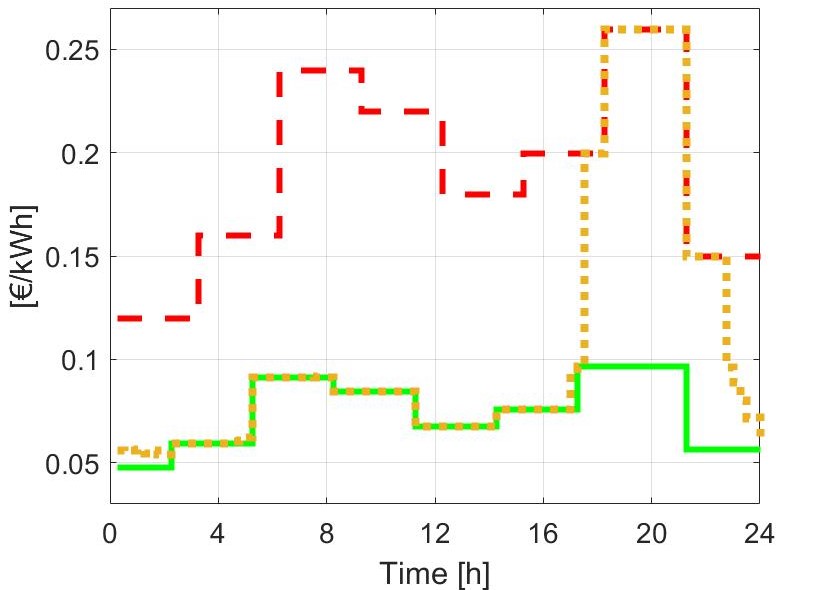}}
\subfloat[]{\includegraphics[width=0.5\linewidth]{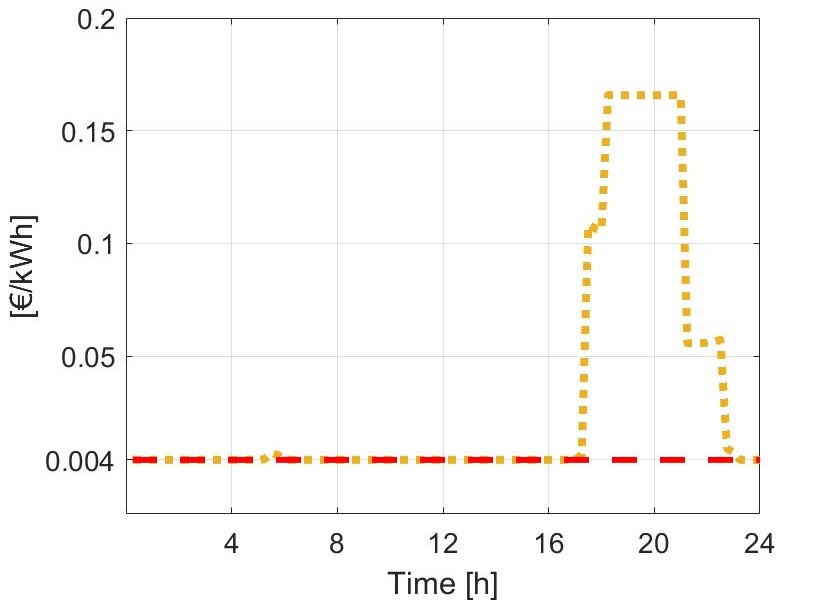}}
\caption{Case study 1: Distributed algorithm results: (a) eAG total active power output; (b) eAG overall active up power reserve; (c) Energy prices trends: buying (dashed), selling (solid) and internal energy price $\boldsymbol{\lambda^p}/\tau$ (dotted) (d) Up reserve energy price trends: main utility (dashed), internal  $\boldsymbol{\lambda^{r_p \uparrow}}\hspace{-1.8mm}/\tau$ (dotted).}
\label{fig:ApoutRup}
\end{figure}
\begin{figure}[!t]
\centering
\subfloat[]{\includegraphics[width=0.5\linewidth]{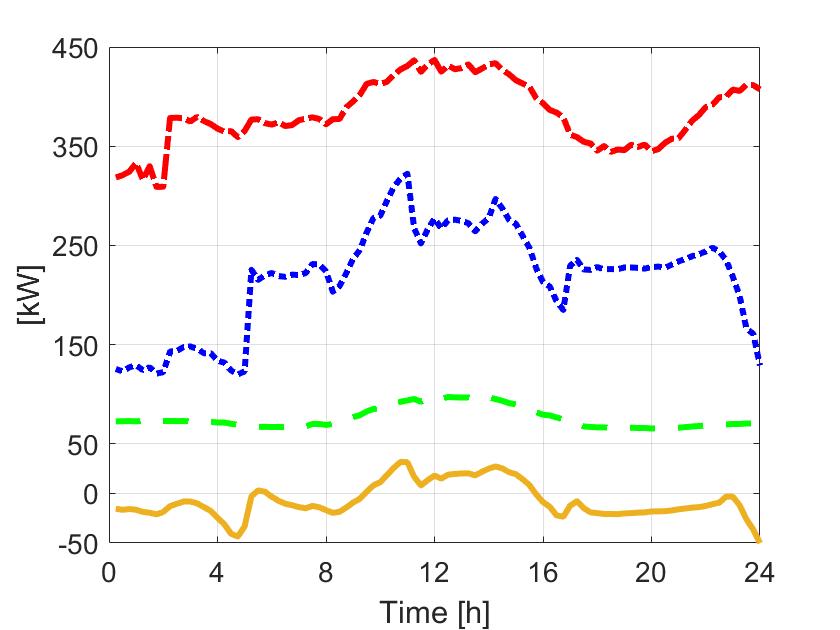}}
\subfloat[]{\includegraphics[width=0.5\linewidth]{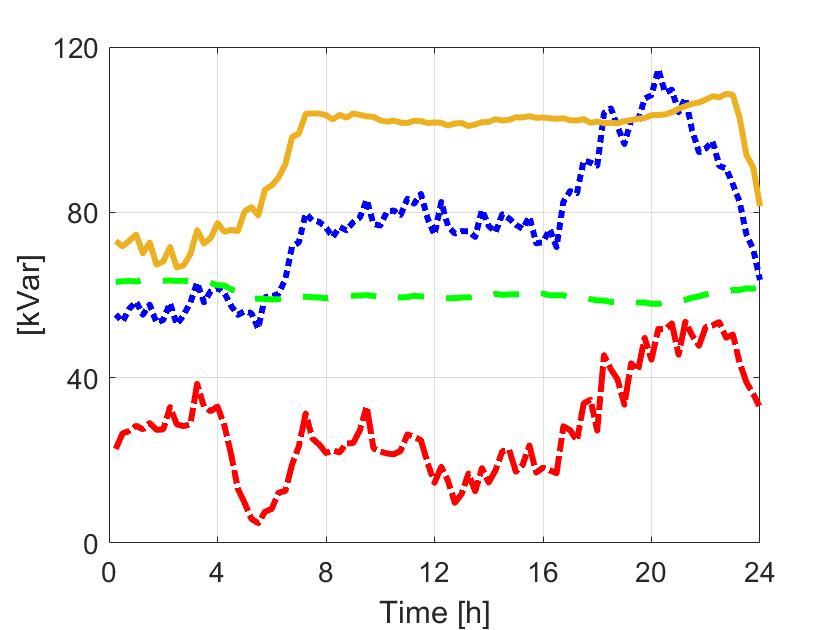}} \\
\subfloat[]{\includegraphics[width=0.49\linewidth]{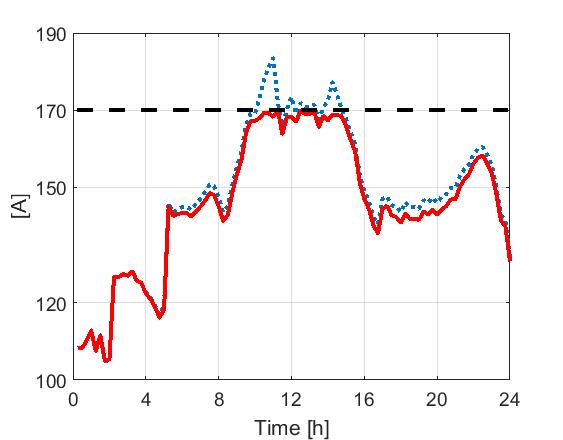}}
\subfloat[]{\includegraphics[width=0.5\linewidth]{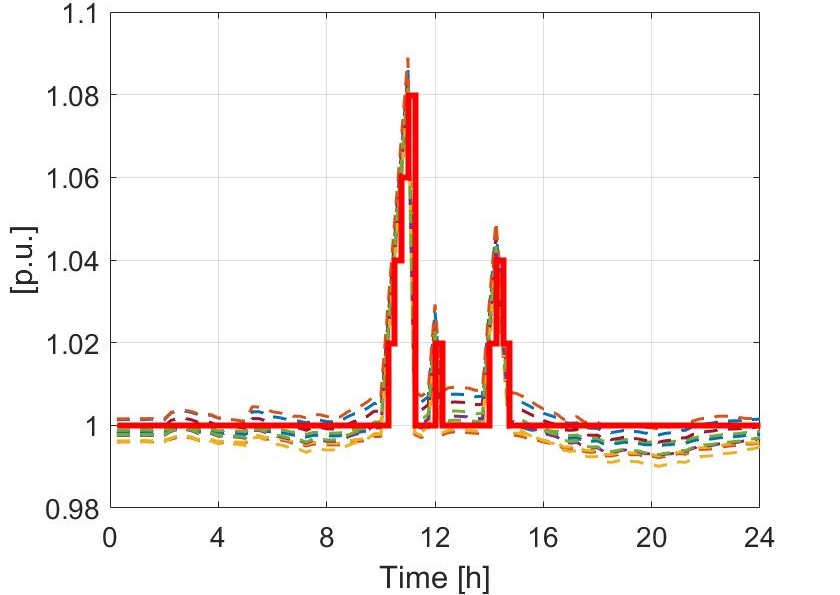}}
\caption{Case study 1: (a)-(b) MGs active and reactive power flows: MG1 (dashed-dotted), MG2 (dotted), MG3 (dashed) and MG4 (solid);
(c) Line 7-8 current magnitude in case current bounds are considered (solid line) and in case they are not (dotted line);  (d) Nodal voltages in case current bounds are considered: $V_1$ (solid),  $V_{2,\hdots,13}$ (dashed).}
\label{fig:VIPQ}
\end{figure}
Once the active power profile of each MG has been scheduled, Phase 2 can be performed to check if electrical feasibility is compromised. The corresponding optimization problem, \eqref{centr_opf}, has been solved using the BMIBMB solver in MATLAB, which can deal with non-convex problems with mixed-integer variables. The current and voltage trends have been simulated through the MATPOWER environment, a well-known power flow simulation tool \cite{zimmerman2011matpower}.
The final active and reactive power trends for each MG, scheduled in Phase 1 and Phase 2 respectively, are shown in Figure \ref{fig:VIPQ}(a)-(b). In the considered case study, it was not necessary to vary the active power outputs scheduled in the Phase 1 since the electrical feasibility was not compromised. \\
Considering again Figure~\ref{fig:ieee13}, it can be noted that line 7-8 is one of the most critical ones since MG1 and MG2 inject power through it, being also the MGs with larger production capabilities. In Figure \ref{fig:VIPQ}(c) the current magnitude profile in line 7-8 is depicted both in case the current constraints are included in the  formulation of Phase 2 and in case they are not. It is evident that the scheduled power flows would violate the maximum current bound during the central part of the day if line current limits are not considered in the MGs management. On the contrary, if the currents constraints are included in Phase 2, the AGS manages to ensure the electrical feasibility and this is achieved by raising up the OLTC slack voltage, see Figure \ref{fig:VIPQ}(d). The OLTC manipulation therefore allows to make the current evolve inside the predefined range, without the necessity of decreasing the already scheduled active power flows of MG1 and MG2.\\
Finally, the generation units power set-points are defined according to Phase 3(a), since active power variations have not been requested in Phase 2.  Considering for brevity only MG2, its final active and reactive power set-points of MG2 are shown in Figure \ref{fig:REsch}.
\\
\subsection*{Analysis of convergence and sensitivity of the distributed algorithm of Phase 1}
\medskip
\textbf{Algorithm 1} has been tested also considering larger numbers of MGs in order to assess the scalability of the proposed approach. Two additional cases have been considered: one with 8 MGs and one with 12 MGs. For the sake of clarity, the structure of the MGs described in Table \ref{tab:MGs_gen} have been reused. Precisely, for in 8-MGs case, in addition to the previously defined MGs, also MG5...MG8 are considered and placed at nodes 8, 7, 5 and 3, respectively (see again Figure \ref{fig:ieee13}). Moreover, MG5...MG8 are defined equal to MG1...MG4, respectively, but with the initial value of the SOCs of all batteries at 80\%. Since the number of MGs is doubled with respect to the previous case, the minimum amount of required active reserve has been set to 200 kW.
 In the 12-MGs case, the MG9...MG12 are also introduced, and placed at nodes 11, 2, 12 and 5, respectively. MG9...MG12 have been defined with the same generations units of MG1...MG4, respectively, and with initial value of the SOCs of all batteries at 30 \%. Moreover, in this case the aggregation must satisfy the minimum active reserve requirement of 300 kW.
 \\In Figure \ref{fig:conv} the differences between the optimal value of the cost function of the centralized solution and the one of the distributed algorithm over the number of iterations are shown, considering the 4-MGs, 8-MGs and 12-MGs test cases. It is evident that the distributed algorithm achieves the same optimal objective of the centralized solution in about 100 iterations, even though a larger numbers of MGs are considered. This result has been obtained by tuning the step-size update of \textbf{Algorithm 1} based on the aggregation-size. This has been set to $\mu=1e^{-5}$ in the 4-MGs case, $\mu=0.95e^{-5}$ in the 8-MGs case and $\mu=0.7e^{-5}$  in the 12-MGs case. Additional guidelines on the choice of the ADMM step-size can be adopted, as explained in \cite{ghadimi2015optimal}.\\ The test-cases have been simulated using a laptop with an Intel Core i7-6500u processor and 8 GB of ram. The average computational time for computing the first step of \textbf{Algorithm 1} was $t_{1}=3 \,s$, while it was $t_{2}=5 \,s$ for the second and third steps, resulting in a total computational time of $t_{k}=8 \,s$ for each iteration of \textbf{Algorithm 1}. Therefore, considering that about 100 iterations were needed to reach convergence, the total computational time for the algorithm to converge was $t_{tot}=800\, s$ to schedule the whole daily power production and reserve provision. 
These results witness the potentiality of the approach as the total computational time to obtain the optimal schedule does not increase with the number of MGs, even though the size of the problem considerably grows up. Indeed, considering both the MGs and the AGS optimization problems,  3456 optimization variables must be considered for the 4-MGs case, 6624 for the 8-MGs case and 9792 for the 12-MGs case.

\begin{figure}[!t]
\centering
{\includegraphics[width=0.7\linewidth]{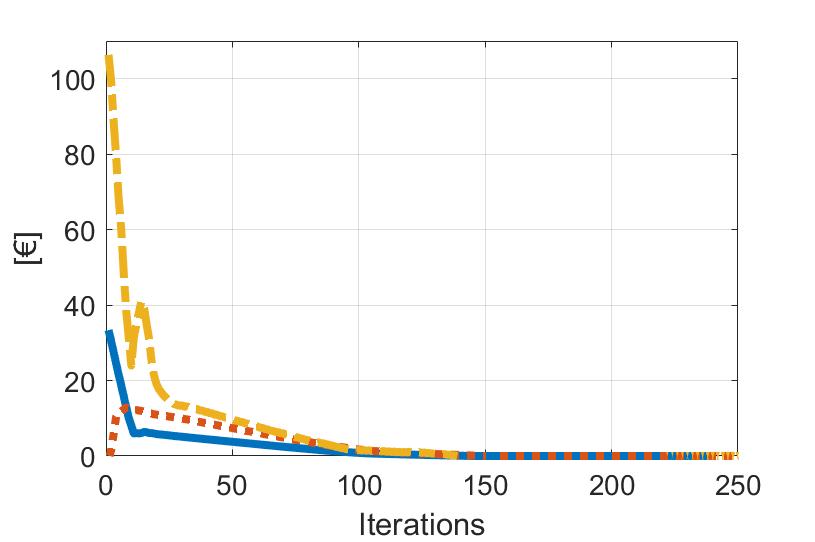}}
\caption{Case study 1: Optimal objective gap between the centralized and distributed solution of Phase 1,$(\mathbf{J^{\textbf{\tiny{AG}},*}} - \mathbf{J^{\textbf{\tiny{AG}},k}})$:  4 MGs (solid line), 8 MGs (dotted line), 12 MGs (dashed line).}\label{fig:conv}
\end{figure}

\subsection{Case study 2: IEEE 37-bus system}\label{sec:52}
\begin{figure}[t!]
	\centering
	\includegraphics[width=0.9\linewidth]{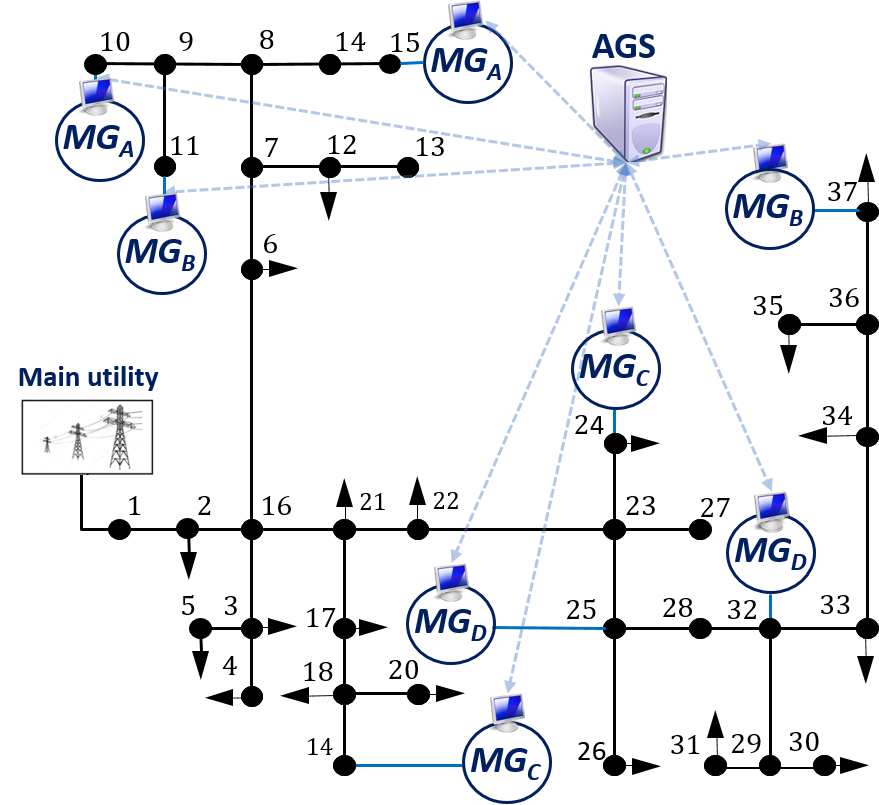}
	\caption{Case study 2: eAG network topology: IEEE 37 bus system}
	\label{fig:ieee37}
\end{figure}

As shown in Figure~\ref{fig:ieee37}, the eAG is composed of 8 MGs and 17 loads. Four different MGs structures are considered and their characteristics are reported in Table \ref{tab:MGs_gen}. All the storage systems are supposed to start with a random initial  SOC, limited between $20\%$ and $80\%$. The loads active power trends are shown in Figure \ref{fig:forecasts}(a); concerning the reactive power trends, it is assumed that all the loads have constant power factor $0.8$. In Figure \ref{fig:forecasts}(b), the day-ahead energy selling and buying prices are shown, while, for simplicity, the up and down power reserve prices are kept equal to a constant value for the whole day, i.e.  $\rho^{\uparrow}_{r_p}=\rho^{\downarrow}_{r_p}= 0.004 $ \euro /kWh. Concerning the ancillary services provision, the whole eAG must provide at least 200 kW of up and down power reserve at each instant, while the line currents are limited by $160$A and the nodal voltages by a maximum variation of $10\%$.
The method has been implemented in MATLAB, using the solvers CPLEX for \emph{Phases 1} and \emph{3} and \emph{fmincon} for \emph{Phase 2}.
The overall optimization problem has been solved also in a centralized fashion for comparison. In the considered case, ADMM converges to the global optimum in about 150 iterations.
Figures \ref{fig:ApoutRup}(a)-(b) show the scheduled eAG output power and the provided up-reserve trends achieved using both the centralized and the distributed approaches: it is evident that the same optimal solution is obtained and that the minimum required amount of power reserve is always globally provided. This means that the AGS is able to optimally manage the eAG without any information about the MGs composition and internal optimization problems. Once the active power profile of each MG has been scheduled, \emph{Phase 2} has been performed to check if electrical feasibility is compromised by the optimal active power flows computed in \emph{Phase 1}.
\begin{table}
	\caption {Case study 2: MGs generation units}
	\label{tab:MGs_gen}
	\centering
	\begin{tabular}{p{1.2cm} p{1.05cm} p{1.5cm} p{1.1cm} p{1.6cm} }
		\toprule
		&  & & &\\[-1em]
		Owner & $(\underline{p},\,\bar{p})$ & $(\underline{q},\,\bar{q})$ & $(C^b,\bar{E}^b)$ & Costs  \\ 		\midrule
		MG$_\text{A}$(Gen.) & {$(20, 300)$} & $\pm\, \sqrt{(330 \text{\,-\,} p^g)}$ & $-$ & $(0.02, 6, 30)e^{-3}$ \\           MG$_\text{A}$(Gen.) & $(20, 300)$ & $ \pm\, \sqrt{(330 \text{\,-\,} p^g)}$ & $-$ & $(0.03, 6, 30)e^{-3}$ \\
		MG$_\text{A}$(Res.) & $(0, 80) $& $\pm $80  & $-$ & $-$ \vspace{1mm} \\
		MG$_\text{B}$(Gen.) & $(20, 300)$ & $\pm\, \sqrt{(330 \text{\,-\,} p^g)}$ & $-$ & $(0.04, 1, 30)e^{-3}$ \\
		MG$_\text{B}$(Batt.) & $\pm$ 30 & $\pm$ 30 &$ (40,20e^3)$ & $(0.95e^{-4}, \,1)$  \\
		MG$_\text{B}$(Res.) & $(0, 150)$ & $\pm 100 $ & $-$ & $-$ \vspace{1mm}\\
		MG$_\text{C}$(Gen.) & $(10, 80)$ & $ \pm\, \sqrt{(90 \text{\,-\,} p^g)}$ & $-$ & $(0.01, 3, 30)e^{-3}$ \\
		MG$_\text{C}$(Res.) & $(0, 30) $ & $\pm$ 30  & $-$ & $-$ \vspace{1mm}\\
		MG$_\text{D}$(Batt.) & $\pm$ 40 & $\pm$ 40 & $(50,\,20e^3)$ & $(1.15e^{-4}, 1)$ \\
		MG$_\text{D}$(Batt.) & $\pm$ 30 & $\pm$ 30 & $(40,\,10e^3) $ & (1.10$e^{-4}$, 2) \\
		MG$_\text{D}$(Res.) & (0, 60) & $\pm$ 50  & $-$ & $-$ \\
		\bottomrule
	\end{tabular}
	\setlength\belowcaptionskip{-40pt}
\end{table}
\begin{figure}[!t]
	\vspace{-7mm}
	\centering
	\subfloat[]{\includegraphics[width=0.5\linewidth]{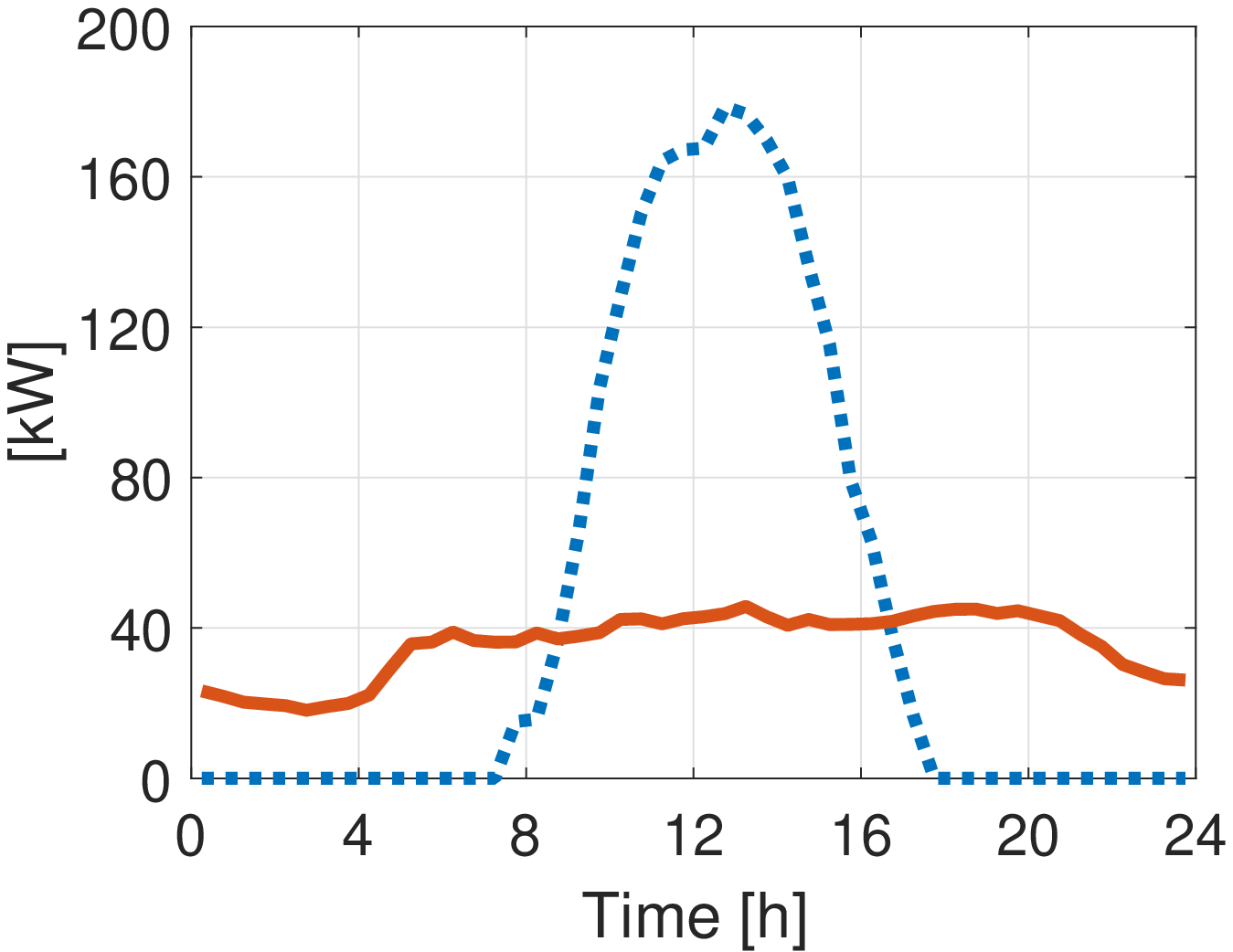}}
	\subfloat[]{\includegraphics[width=0.5\linewidth]{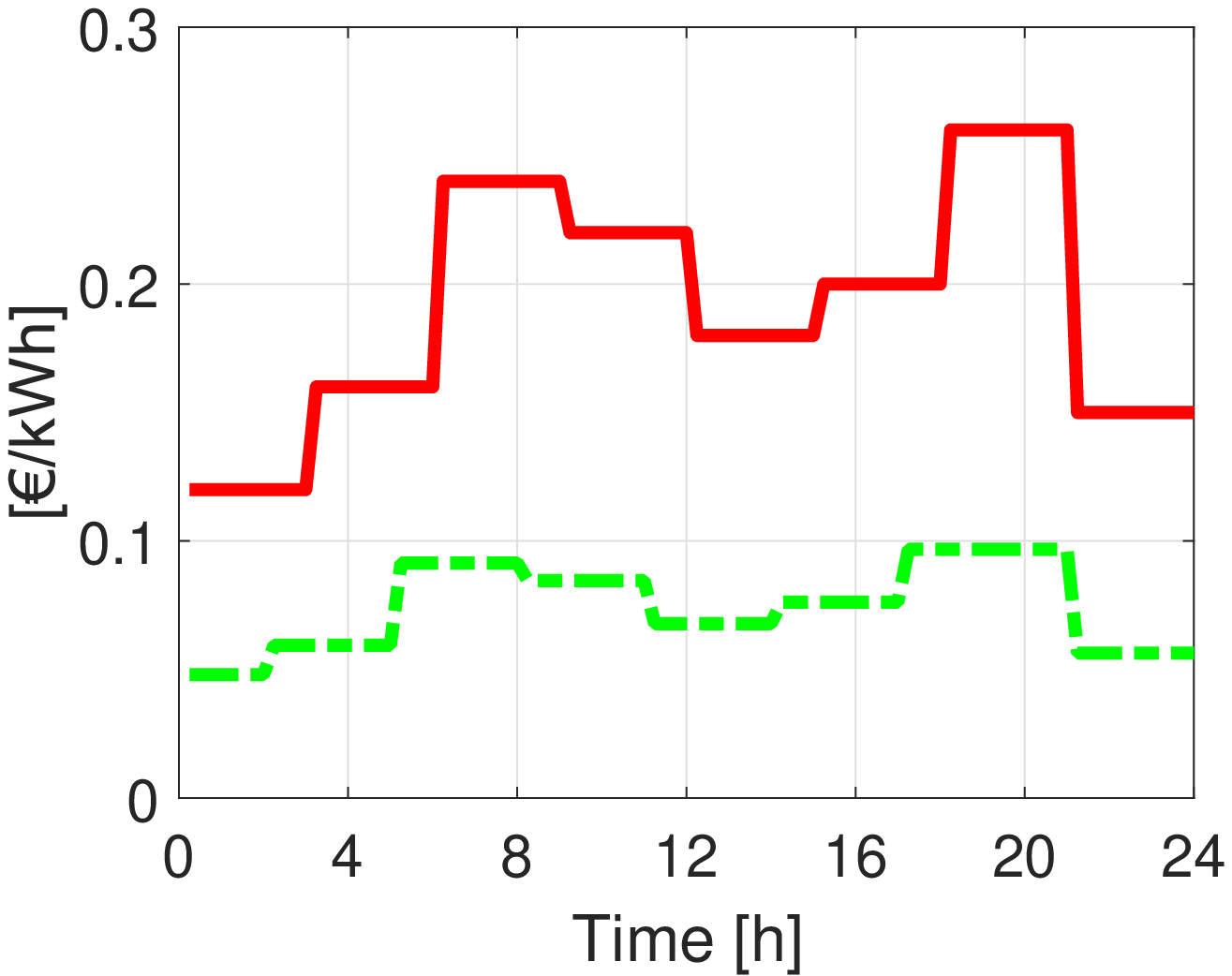}}
	\setlength\belowcaptionskip{-14pt}
	\caption{Case study 2: (a) eAG loads active power demand; (b) Day-ahead energy prices: buying price (solid red), selling price (dashed green).
	}
	\label{fig:forecasts}
\end{figure}
\begin{figure}[!t]
	\centering
	\subfloat[]{\includegraphics[width=0.5\linewidth]{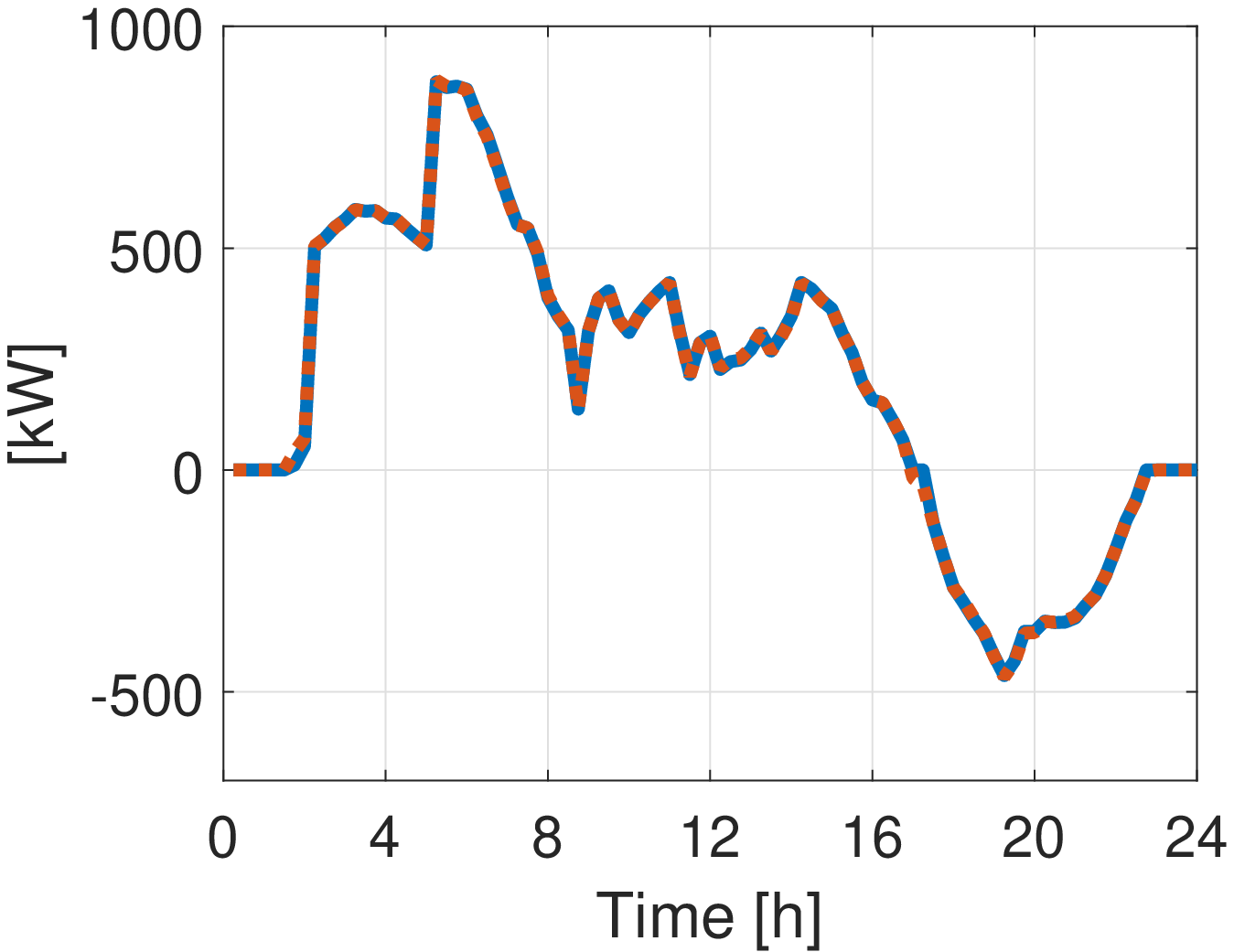}}
	\subfloat[]{\includegraphics[width=0.5\linewidth]{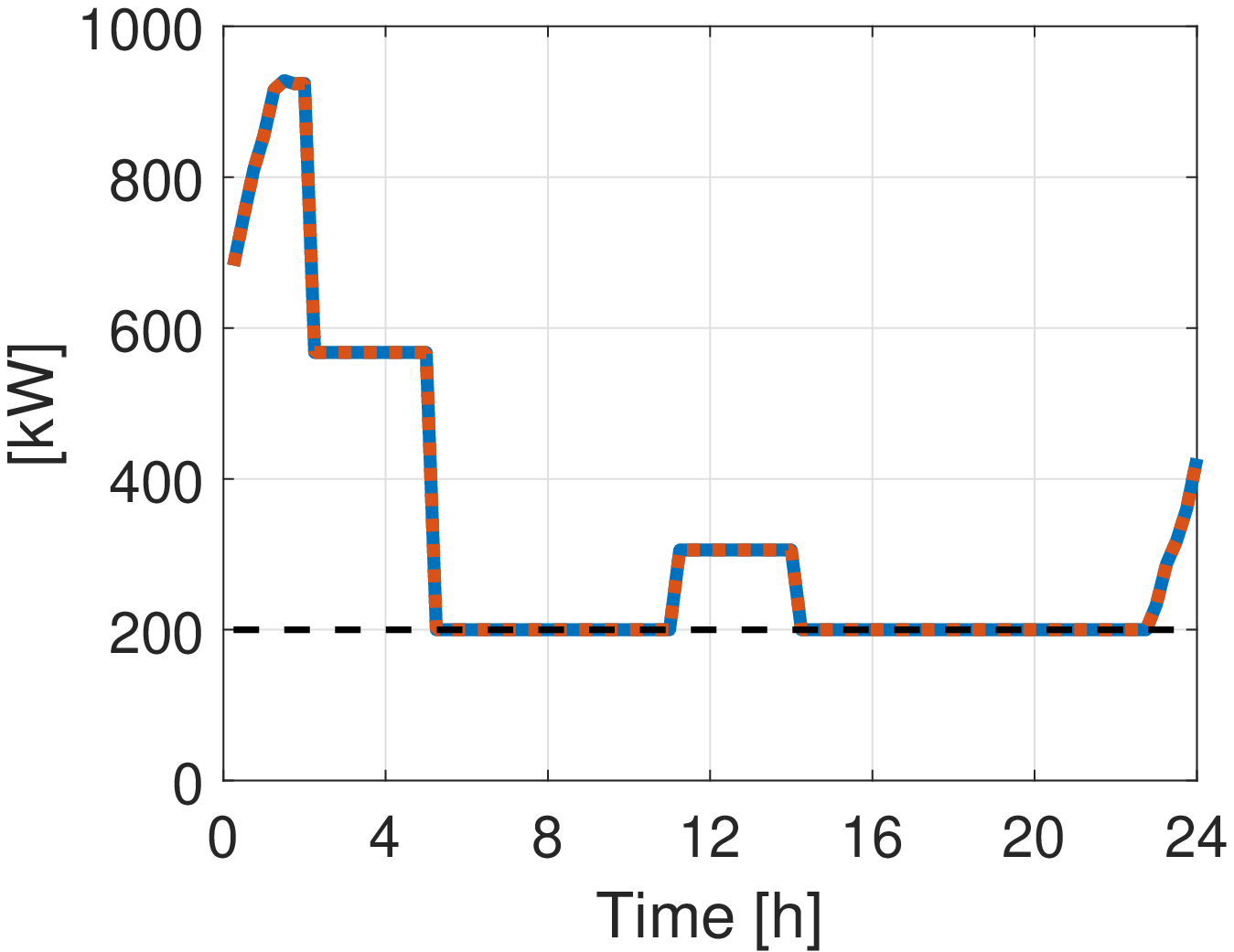}}
	\\[-3mm]
	\subfloat[]{\includegraphics[width=0.5\linewidth]{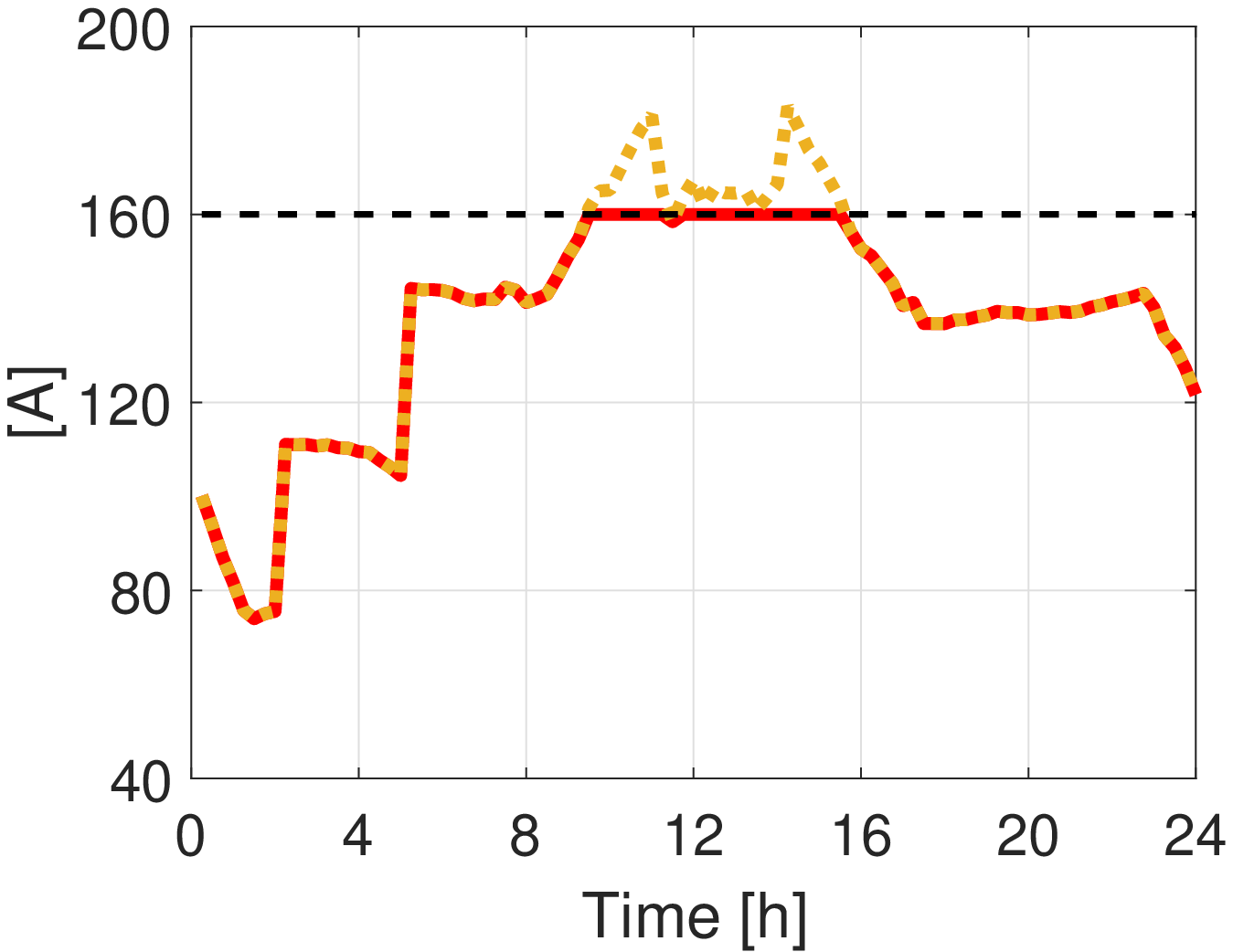}}
	\subfloat[]{\includegraphics[width=0.5\linewidth]{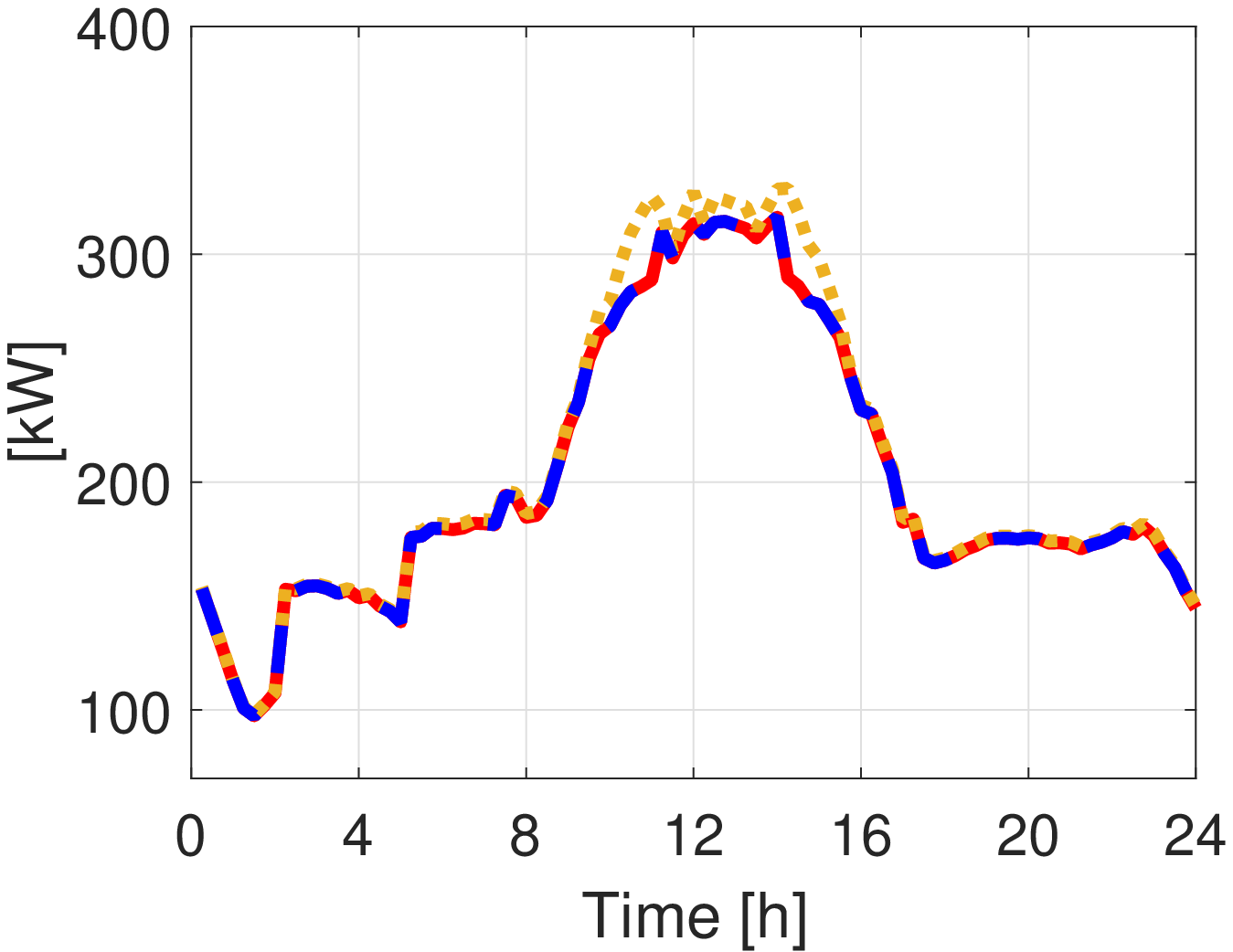}}
	\setlength\belowcaptionskip{-12pt}
	\caption{Case study 2: Centralized (dashed blue) vs distributed solution (solid red): (a) eAG active power output, (b) eAG active up power reserve; (c) Line 7-8 current magnitude in case constraint \eqref{I_constr} is considered (solid red) or  not (dotted yellow); (d) Active power of MG$_{\text{B}}$ attached to node 11 obtained by \emph{Phase 1} (dotted yellow), \emph{Phase 2} (dashed red)  and \emph{Phase 3} (solid blue).
	}
	\label{fig:ApoutRup}
\end{figure}

\begin{figure}[!t]
	\centering
	\subfloat[]{\includegraphics[width=0.5\linewidth]{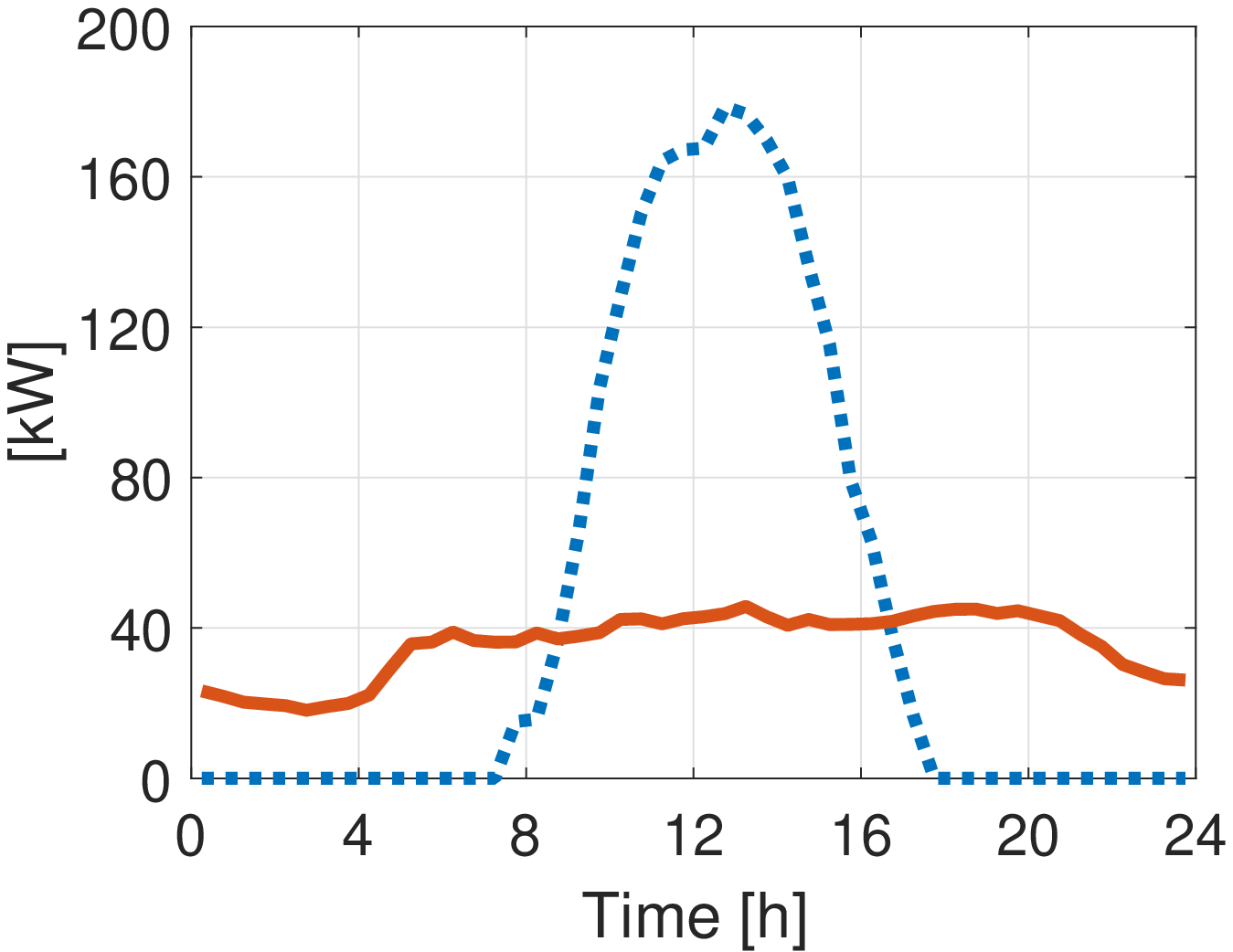}}
	\subfloat[]{\includegraphics[width=0.5\linewidth]{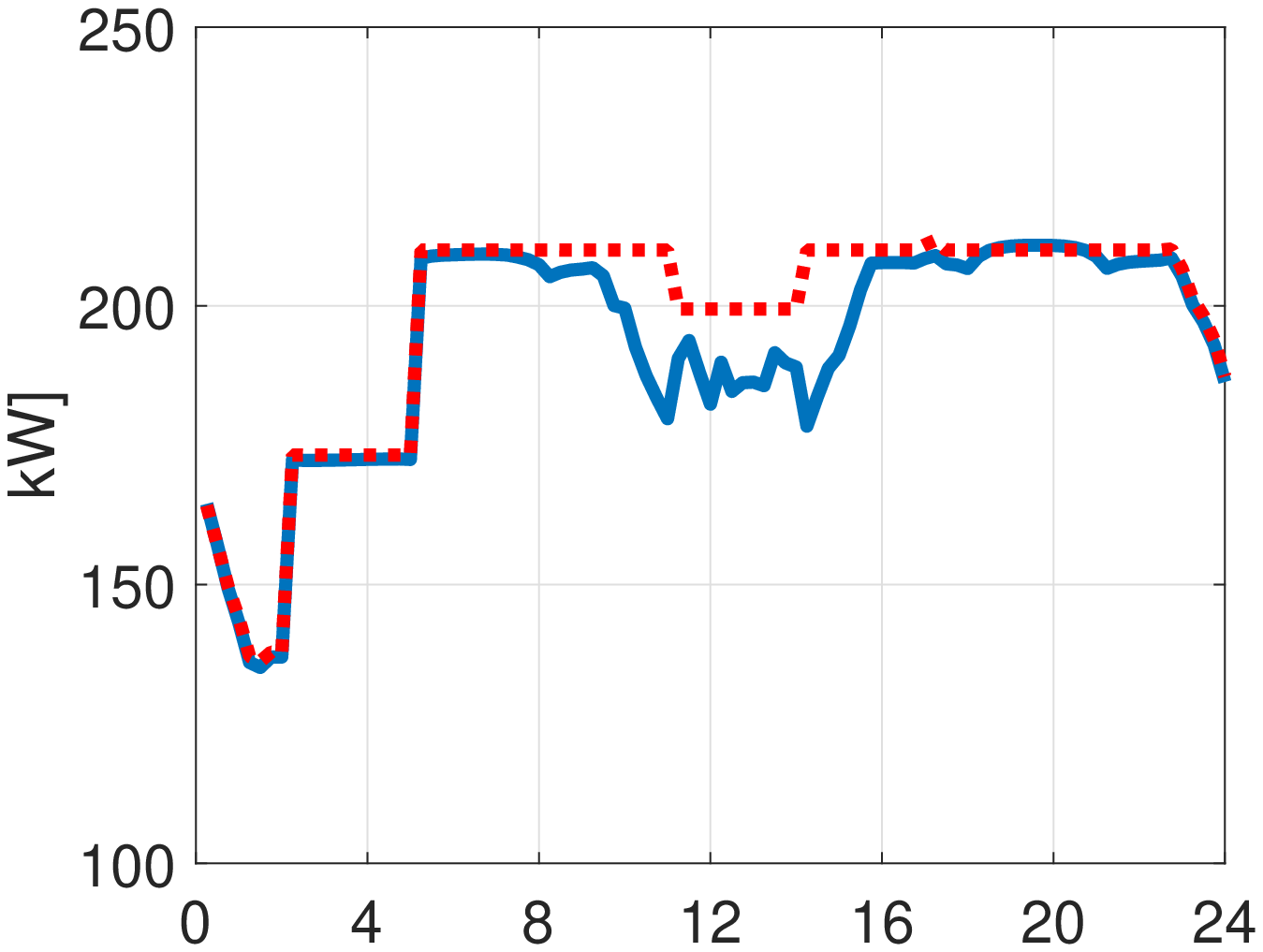}}
	\\[-3mm]
	\subfloat[]{\includegraphics[width=0.5\linewidth]{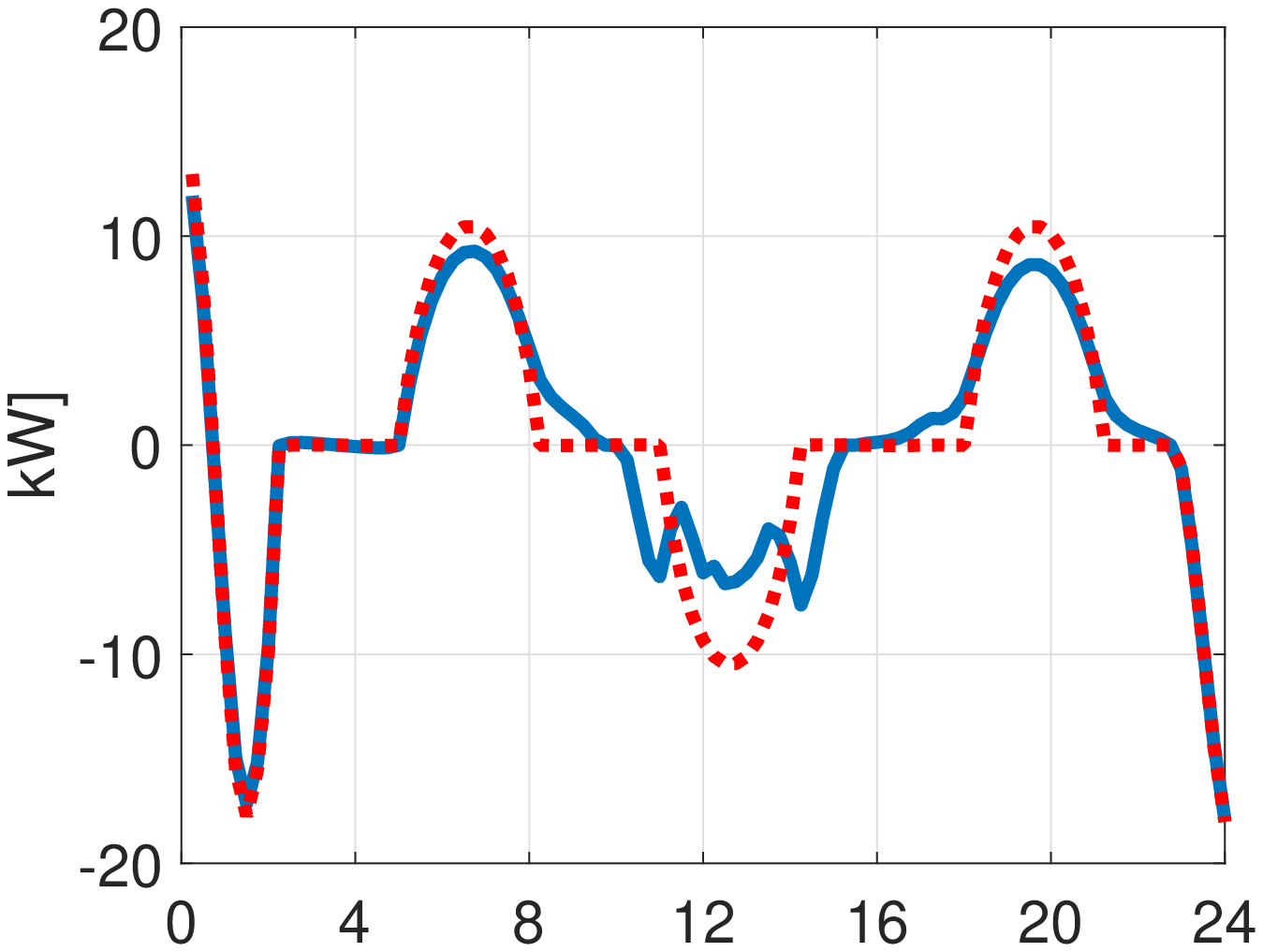}}
	\subfloat[]{\includegraphics[width=0.5\linewidth]{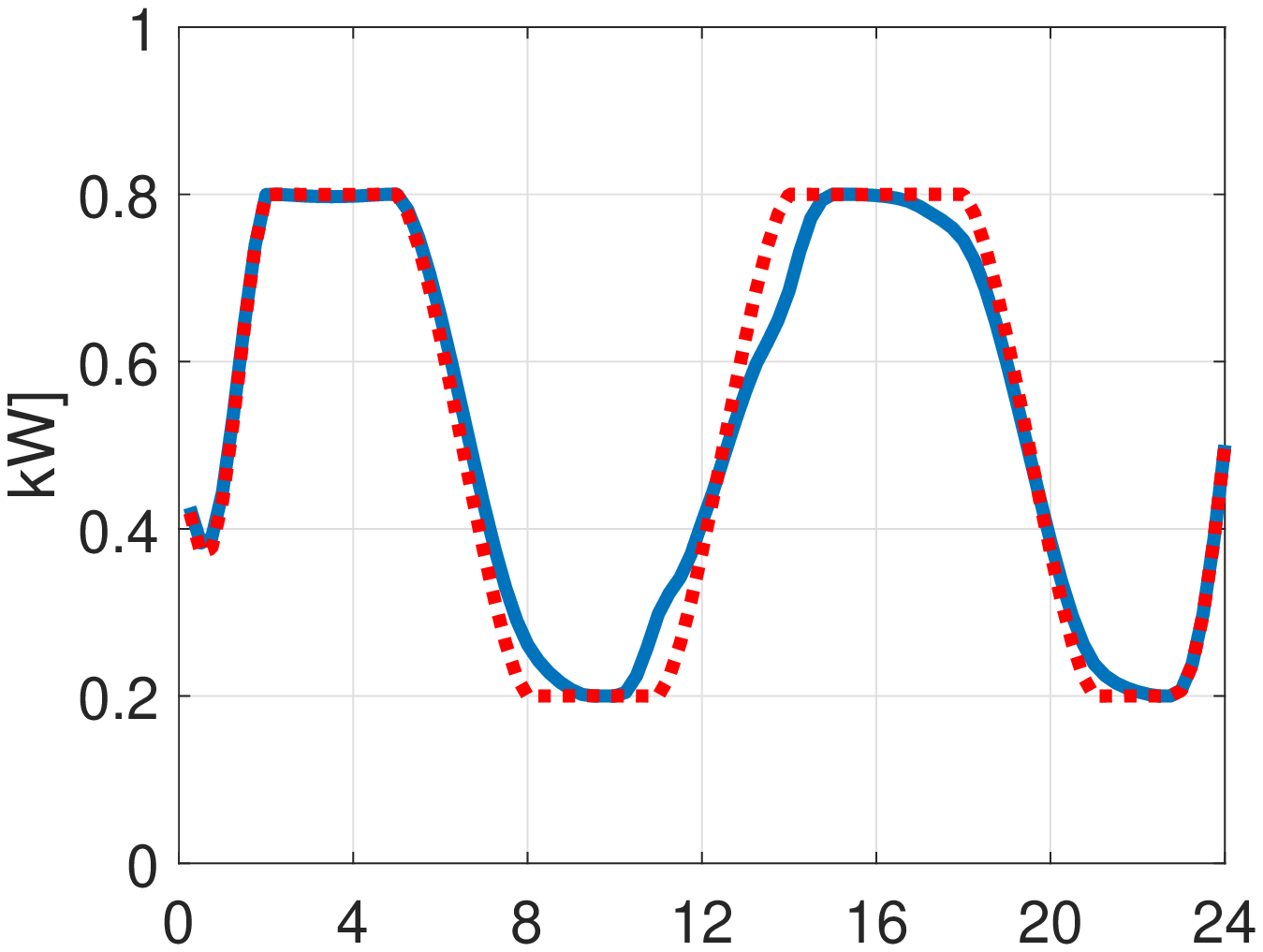}}
	\setlength\belowcaptionskip{-12pt}
	\caption{Case study 2: MG$_{\text{B}}$: (a) Load absorption (solid red) and PV generation  (dashed blue), (b) Generator active power output by Phase 1 (dashed red) and by Phase 3 (solid blue); (c) Battery active power output by Phase 1 (dashed red) and by Phase 3 (solid blue); (d) Battery SOC by Phase 1 (dashed red) and by Phase 3 (solid blue).
	}
	\label{fig:MG11_Phase123}
\end{figure}

%
Considering Figure~\ref{fig:ieee37}, it can be noted that line 7-8 is one of the most critical ones since three MGs  inject power through it. In Figure \ref{fig:ApoutRup}(c) the current magnitude profile in line 7-8 is shown both in case the current constraints are considered and in case they are not. It is evident that the scheduled MGs power flows violate the maximum current bound during the central part of the day if line current limits are not considered. On the contrary, if \emph{Phase 2} is performed, the electrical feasibility can be achieved by asking the MGs a small active power reduction. For the reported numerical results, the MGs placed at nodes 10, 11 and 15 are asked to reduce their power; in Figure \ref{fig:ApoutRup}(d) the output power of the MG attached to node 11 is shown both after \emph{Phase 1} (dotted) and after \emph{Phase 2} (dashed). At this stage, \emph{Phase 3} is executed such that the MGs reschedule their units according to the active and reactive power set-points from \emph{Phase 1} and \emph{Phase 2}. The MGs are able to track the power references and provide the minimum requested power reserve given the fact that very small active power variations are expected to be requested. In Figure \ref{fig:ApoutRup}(d) the final output power of the MG attached to node 11 is shown, which completely overlaps the profile obtained after \emph{Phase 2}. Finally, Figure \ref{fig:MG11_Phase123} reports the profiles of the internal units of the MG attached to node 11. It is worth noticing that as in Phase 2 it was requested to the MG to decrease the output power to avoid over-current issues, the MG then in Phase 3 reschedules its units minimizing its production costs to maintain the output power references received by Phase 1 and Phase 2. As stated in the paper, Phase 3 however ensures that the reserve constraints are respected by the whole eAG.\\
The scalability of {Algorithm 1} has been tested  considering two additional cases: one with 12 MGs and one with 16 MGs, connected to the IEEE 37-bus system. Figure \ref{fig:conv} shows the differences between the optimal values of the cost function of the centralized solution and the one of the distributed algorithm for the 8-MGs, 12-MGs and 16-MGs test cases. The distributed algorithm achieves the same optimal objective of the centralized solution in about 150 iterations, even though a larger number of MGs are considered. This result has been obtained by properly tuning the parameter $\mu$ of {Algorithm 1} based on the size of the aggregation. Additional guidelines on the choice of the ADMM step-size can be found in \cite{ghadimi2015optimal}.
The average computational time for computing the first step of \textbf{Algorithm 1} was $t_{1}=3 \,s$, while it was $t_{2}=5 \,s$ for the second and third steps, resulting in a total computational time of $t_{k}=8 \,s$ for each iteration of \textbf{Algorithm 1}. Therefore, considering that about 150 iterations were needed to reach convergence, the total computational time for the algorithm to converge was $t_{tot}=1200\, s$ to schedule the whole daily power production and reserve provision. 
These results witness the potentiality of the approach as the total computational time to obtain the optimal schedule does not increase with the number of MGs, even though the size of the problem considerably grows up. Indeed, considering both the MGs and the AGS optimization problems, 6624 for the 8-MGs case, 9792 for the 12-MGs case and 13248 for the 16-MGs case.

\begin{figure}[!t]
	\centering
	{\includegraphics[width=0.7\linewidth]{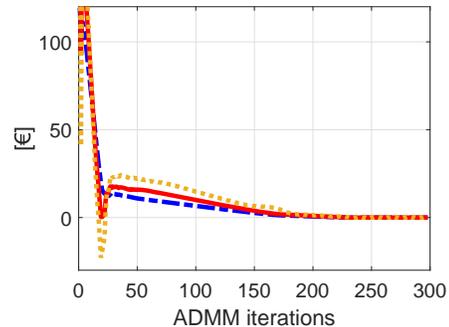}}
	\setlength\belowcaptionskip{-10pt}
	\caption{Case study 2: Optimal objective gap between the centralized and distributed solution of \emph{Phase 1},$(\mathbf{J^{\textbf{\tiny{AG}},*}} - \mathbf{J^{\textbf{\tiny{AG}},k}})$:  8 MGs (dashed blue), 12 MGs (solid red), 16 MGs (dotted yellow).}\label{fig:conv}
\end{figure}

\vspace{-1mm}
\section{Conclusion}\label{sec:6}
In this paper, an approach for the day-ahead  scheduling and control of microgrids grouped to form a electrical aggregation has been described. The goal is to allow microgrids to provide ancillary services to the main utility, e.g. in terms of power reserve for frequency control. The method has been structured in three phases. Scalability of the optimization problem and confidentiality requirements have been guaranteed by resorting to a suitable distributed optimization algorithm.\\
The developed solution can be easily extended in many directions. Among them, we recall the possibility to re-initialize and run the algorithm during the day when unexpected conditions occur in term of produced power due to significant weather forecast variations or load changes. A further research direction can concern the use of these ideas to define aggregators of commercial buildings. Future developments will also regard the management of non-convex MGs cost functions, e.g. due to the presence of boolean variables describing the switch on/off of generators and batteries different charge/discharge conditions.
%
\section*{Acknowledgment}
The work of Alessio La Bella and Carlo Sandroni has been financed by the Research Fund for the
Italian Electrical System under the Contract Agreement
between RSE S.p.A. and the Ministry of Economic
Development - General Directorate for Nuclear Energy, Renewable Energy and Energy Efficiency - in compliance
with the Decree of March 8, 2006.

\appendices
\section{Proof of Proposition 1} \label{Prop1_proof}
\begin{proof}
Proposition 1 holds $\; \forall t \, \in \, \{1, \hdots,N\}$ and therefore the time index is neglected in the following.\\ The function
\begin{align*}
g^p_z(z^p)= &- \,{\rho_e^{s}}\;{\max}\,({p}^{{{\text{\tiny{AG}}}}}\tau ,{0}) +{\rho_e^{b}}\;{\max}(-p^{{{\text{\tiny{AG}}}}}\tau ,{0}) \,+\\&-{\rho_{r_p}^{\uparrow}}\;{r}_p^{\text{\tiny{AG}}\uparrow} \tau  - {\rho_{r_p}^{\downarrow}}  \; {r}_p^{\text{\tiny{AG}}\downarrow}\tau
\end{align*}
is said to be convex over $\mathcal{Z}^p$ if and only if
\begin{align}
g^p_z(z_2^p) \, \geq \, g^p_z(z_1^p) \, + \, (z_2^p-z_1^p)' \,\nabla g^p_z(z_1) \quad \forall z_1^p,z_2^p \in \mathcal{Z}^p
\label{convexity}
\end{align}
Moreover, if the inequality \eqref{convexity} is strict with $z^p_1 \neq z^p_2$, then $g_z(z^p)$ is said to be strictly convex over $\mathcal{Z}^p$.\\
Take two instances $z_1^p=[({p}^{{{\text{\tiny{AG}}}}})_1,({r}_{p}^{\text{\tiny{AG}}\uparrow})_1,({r}_{p}^{\text{\tiny{AG}}\downarrow})_1]'$ and $z^p_2=[({p}^{{{\text{\tiny{AG}}}}})_2,({r}_{p}^{\text{\tiny{AG}}\uparrow})_2,({r}_{p}^{\text{\tiny{AG}}\downarrow})_2]'$.
If $({p}^{{{\text{\tiny{AG}}}}})_2\,\leq\,0\,\leq\,({p}^{{{\text{\tiny{AG}}}}})_1$,  it can be verified that \eqref{convexity} becomes $(\rho_e^b - \rho^s_e) \cdot\,({p}^{{{\text{\tiny{AG}}}}})_2 \leq 0$. Since $({p}^{{{\text{\tiny{AG}}}}})_2 \leq 0$, this implies that $g_z(z^p)$ is convex if and only if $\rho_e^b \geq \rho_e^s$. The same result  can be easily obtained in case $({p}^{{{\text{\tiny{AG}}}}})_2\,\leq\,0\,\leq\,({p}^{{{\text{\tiny{AG}}}}})_1$. On the other hand, if $({p}^{{{\text{\tiny{AG}}}}})_1$ and $({p}^{{{\text{\tiny{AG}}}}})_2$ have the same sign, e.g. $ ({p}^{{{\text{\tiny{AG}}}}})_1\,\geq\,({p}^{{{\text{\tiny{AG}}}}})_2 \, \geq 0$, \eqref{convexity} collapses in a trivial equality. This implies that $g_z(z^p)$ is convex but not strictly.
\end{proof}
\vspace{-5mm}
\bibliographystyle{IEEEtran}
\bibliography{IEEEabrv,mybiblio}

\begin{thebibliography}{10}
\providecommand{\url}[1]{#1}
\csname url@samestyle\endcsname
\providecommand{\newblock}{\relax}
\providecommand{\bibinfo}[2]{#2}
\providecommand{\BIBentrySTDinterwordspacing}{\spaceskip=0pt\relax}
\providecommand{\BIBentryALTinterwordstretchfactor}{4}
\providecommand{\BIBentryALTinterwordspacing}{\spaceskip=\fontdimen2\font plus
\BIBentryALTinterwordstretchfactor\fontdimen3\font minus
  \fontdimen4\font\relax}
\providecommand{\BIBforeignlanguage}[2]{{%
\expandafter\ifx\csname l@#1\endcsname\relax
\typeout{** WARNING: IEEEtran.bst: No hyphenation pattern has been}%
\typeout{** loaded for the language `#1'. Using the pattern for}%
\typeout{** the default language instead.}%
\else
\language=\csname l@#1\endcsname
\fi
#2}}
\providecommand{\BIBdecl}{\relax}
\BIBdecl

\bibitem{lasseter2002microgrids}
R.~H. Lasseter, ``Microgrids,'' in \emph{Power Engineering Society Winter
  Meeting, 2002. IEEE}, vol.~1.\hskip 1em plus 0.5em minus 0.4em\relax IEEE,
  2002, pp. 305--308.

\bibitem{guerrero2011hierarchical}
J.~M. Guerrero, J.~C. Vasquez, J.~Matas, L.~G. De~Vicu{\~n}a, and M.~Castilla,
  ``Hierarchical control of droop-controlled {AC} and {DC} microgrids-a general
  approach toward standardization,'' \emph{IEEE Transactions on Industrial
  Electronics}, vol.~58, no.~1, pp. 158--172, 2011.

\bibitem{la2017hierarchical}
A.~La~Bella, S.~Raimondi~Cominesi, C.~Sandroni, and R.~Scattolini,
  ``Hierarchical predictive control of microgrids in islanded operation,''
  \emph{IEEE Transactions on Automation Science and Engineering}, vol.~14,
  no.~2, pp. 536--546, 2017.

\bibitem{parisio2014model}
A.~Parisio, E.~Rikos, and L.~Glielmo, ``A model predictive control approach to
  microgrid operation optimization,'' \emph{IEEE Transactions on Control
  Systems Technology}, vol.~22, no.~5, pp. 1813--1827, 2014.

\bibitem{cominesi2016multi}
S.~Raimondi~Cominesi, A.~La~Bella, M.~Farina, and R.~Scattolini, ``A
  multi-layer control scheme for microgrid energy management,'' \emph{IFAC
  Workshop on Control of Transmission and Distribution Smart Grids - CTDSG
  2016, Prague. IFAC-PapersOnLine}, vol.~49, no.~27, pp. 256--261, 2016.

\bibitem{majzoobi2017application}
A.~Majzoobi and A.~Khodaei, ``Application of microgrids in providing ancillary
  services to the utility grid,'' \emph{Energy}, vol. 123, pp. 555--563, 2017.

\bibitem{joos2000potential}
G.~Joos, B.~Ooi, D.~McGillis, F.~Galiana, and R.~Marceau, ``The potential of
  distributed generation to provide ancillary services,'' in \emph{Power
  Engineering Society Summer Meeting, 2000. IEEE}, vol.~3.\hskip 1em plus 0.5em
  minus 0.4em\relax IEEE, 2000, pp. 1762--1767.

\bibitem{qureshi2016economic}
F.~A. Qureshi, I.~Lymperopoulos, A.~A. Khatir, and C.~N. Jones, ``Economic
  advantages of office buildings providing ancillary services with intraday
  participation,'' \emph{IEEE Transactions on Smart Grid}, 2016.

\bibitem{yuen2007feasibility}
C.~Yuen and A.~Oudalov, ``The feasibility and profitability of ancillary
  services provision from multi-microgrids,'' in \emph{Power Tech, 2007 IEEE
  Lausanne}.\hskip 1em plus 0.5em minus 0.4em\relax IEEE, 2007, pp. 598--603.

\bibitem{koliou2014demand}
E.~Koliou, C.~Eid, J.~P. Chaves-{\'A}vila, and R.~A. Hakvoort, ``Demand
  response in liberalized electricity markets: Analysis of aggregated load
  participation in the {German} balancing mechanism,'' \emph{Energy}, vol.~71,
  pp. 245--254, 2014.

\bibitem{yuen2011provision}
C.~Yuen, A.~Oudalov, and A.~Timbus, ``The provision of frequency control
  reserves from multiple microgrids,'' \emph{IEEE Transactions on Industrial
  Electronics}, vol.~58, no.~1, pp. 173--183, 2011.

\bibitem{vrettos2016robust}
E.~Vrettos, F.~Oldewurtel, and G.~Andersson, ``Robust energy-constrained
  frequency reserves from aggregations of commercial buildings,'' \emph{IEEE
  Transactions on Power Systems}, vol.~31, no.~6, pp. 4272--4285, 2016.

\bibitem{boyd2011distributed}
S.~Boyd, N.~Parikh, E.~Chu, B.~Peleato, J.~Eckstein \emph{et~al.},
  ``Distributed optimization and statistical learning via the alternating
  direction method of multipliers,'' \emph{Foundations and
  Trends{\textregistered} in Machine Learning}, vol.~3, no.~1, pp. 1--122,
  2011.

\bibitem{balandat2014contract}
M.~Balandat, F.~Oldewurtel, M.~Chen, and C.~Tomlin, ``Contract design for
  frequency regulation by aggregations of commercial buildings,'' in
  \emph{Communication, Control, and Computing (Allerton), 2014 52nd Annual
  Allerton Conference on}.\hskip 1em plus 0.5em minus 0.4em\relax IEEE, 2014,
  pp. 38--45.

\bibitem{mai2015economic}
W.~Mai and C.~Chung, ``Economic {MPC} of aggregating commercial buildings for
  providing flexible power reserve,'' \emph{IEEE Transactions on Power
  Systems}, vol.~30, no.~5, pp. 2685--2694, 2015.

\bibitem{gkatzikis2013role}
L.~Gkatzikis, I.~Koutsopoulos, and T.~Salonidis, ``The role of aggregators in
  smart grid demand response markets,'' \emph{IEEE Journal on Selected Areas in
  Communications}, vol.~31, no.~7, pp. 1247--1257, 2013.

\bibitem{kim2011two}
H.~Kim and M.~Thottan, ``A two-stage market model for microgrid power
  transactions via aggregators,'' \emph{Bell Labs Technical Journal}, vol.~16,
  no.~3, pp. 101--107, 2011.

\bibitem{dagdougui2014decentralized}
H.~Dagdougui and R.~Sacile, ``Decentralized control of the power flows in a
  network of smart microgrids modeled as a team of cooperative agents,''
  \emph{IEEE Transactions on Control Systems Technology}, vol.~22, no.~2, pp.
  510--519, 2014.

\bibitem{brehm2017distributed}
R.~Brehm, H.~Ramezani, and J.~Jouffroy, ``Distributed coordination of
  energy-storage capacities in virtual microgrids,'' in \emph{European Control
  Conference (ECC'18)}, 2018.

\bibitem{rey2018strengthening}
F.~Rey, X.~Zhang, S.~Merkli, V.~Agliati, M.~Kamgarpour, and J.~Lygeros,
  ``Strengthening the group: Aggregated frequency reserve bidding with
  {ADMM},'' \emph{IEEE Transactions on Smart Grid}, 2018.

\bibitem{guo2017case}
J.~Guo, G.~Hug, and O.~K. Tonguz, ``A case for nonconvex distributed
  optimization in large-scale power systems,'' \emph{IEEE Transactions on Power
  Systems}, vol.~32, no.~5, pp. 3842--3851, 2017.

\bibitem{magnusson2015distributed}
S.~Magn{\'u}sson, P.~C. Weeraddana, and C.~Fischione, ``A distributed approach
  for the optimal power-flow problem based on {ADMM} and sequential convex
  approximations,'' \emph{IEEE Transactions on Control of Network Systems},
  vol.~2, no.~3, pp. 238--253, 2015.

\bibitem{molzahn2017survey}
D.~K. Molzahn, F.~D{\"o}rfler, H.~Sandberg, S.~H. Low, S.~Chakrabarti,
  R.~Baldick, and J.~Lavaei, ``A survey of distributed optimization and control
  algorithms for electric power systems,'' \emph{IEEE Transactions on Smart
  Grid}, vol.~8, no.~6, pp. 2941--2962, 2017.

\bibitem{labellaECC2018}
A.~La~Bella, M.~Farina, C.~Sandroni, and R.~Scattolini, ``Microgrids
  aggregation management providing ancillary services,'' in \emph{2018 European
  Control Conference (ECC), {\rm to appear}}.\hskip 1em plus 0.5em minus
  0.4em\relax IEEE, 2018.

\bibitem{sonmez2013estimation}
Y.~S{\"o}nmez, ``Estimation of fuel cost curve parameters for thermal power
  plants using the {ABC} algorithm,'' \emph{Turkish Journal of Electrical
  Engineering \& Computer Sciences}, vol.~21, no. Sup. 1, pp. 1827--1841, 2013.

\bibitem{zhang2013robust}
Y.~Zhang, N.~Gatsis, and G.~B. Giannakis, ``Robust energy management for
  microgrids with high-penetration renewables,'' \emph{IEEE Transactions on
  Sustainable Energy}, vol.~4, no.~4, pp. 944--953, 2013.

\bibitem{cominesi2018two}
S.~Raimondi~Cominesi, M.~Farina, L.~Giulioni, B.~Picasso, and R.~Scattolini,
  ``A two-layer stochastic model predictive control scheme for microgrids,''
  \emph{IEEE Transactions on Control Systems Technology}, vol.~26, no.~1, pp.
  1--13, 2018.

\bibitem{bertsekas1989parallel}
D.~P. Bertsekas and J.~N. Tsitsiklis, \emph{Parallel and distributed
  computation: numerical methods}.\hskip 1em plus 0.5em minus 0.4em\relax
  Prentice hall Englewood Cliffs, NJ, 1989, vol.~23.

\bibitem{farina2015model}
M.~Farina, A.~Guagliardi, F.~Mariani, C.~Sandroni, and R.~Scattolini, ``Model
  predictive control of voltage profiles in {MV} networks with distributed
  generation,'' \emph{Control Engineering Practice}, vol.~34, pp. 18--29, 2015.

\bibitem{casavola2017command}
A.~Casavola, F.~Tedesco, and M.~Vizza, ``Command governor strategies for the
  online management of reactive power in smart grids with distributed
  generation,'' \emph{IEEE Transactions on Automation Science and Engineering},
  vol.~14, no.~2, pp. 449--460, 2017.

\bibitem{adibi1994reactive}
M.~Adibi and D.~Milanicz, ``Reactive capability limitation of synchronous
  machines,'' \emph{IEEE Transactions on Power Systems}, vol.~9, no.~1, pp.
  29--40, 1994.

\bibitem{marconato2002electric}
R.~Marconato, \emph{Electric power systems}.\hskip 1em plus 0.5em minus
  0.4em\relax CEI, 2002.

\bibitem{burer2012non}
S.~Burer and A.~N. Letchford, ``Non-convex mixed-integer nonlinear programming:
  A survey,'' \emph{Surveys in Operations Research and Management Science},
  vol.~17, no.~2, pp. 97--106, 2012.

\bibitem{houska2016augmented}
B.~Houska, J.~Frasch, and M.~Diehl, ``An augmented lagrangian based algorithm
  for distributed nonconvex optimization,'' \emph{SIAM Journal on
  Optimization}, vol.~26, no.~2, pp. 1101--1127, 2016.

\bibitem{kersting2001radial}
W.~H. Kersting, ``Radial distribution test feeders,'' in \emph{Power
  Engineering Society Winter Meeting, 2001. IEEE}, vol.~2.\hskip 1em plus 0.5em
  minus 0.4em\relax IEEE, 2001, pp. 908--912.

\bibitem{zimmerman2011matpower}
R.~D. Zimmerman, C.~E. Murillo-S{\'a}nchez, and R.~J. Thomas, ``Matpower:
  Steady-state operations, planning, and analysis tools for power systems
  research and education,'' \emph{IEEE Transactions on power systems}, vol.~26,
  no.~1, pp. 12--19, 2011.

\bibitem{ghadimi2015optimal}
E.~Ghadimi, A.~Teixeira, I.~Shames, and M.~Johansson, ``Optimal parameter
  selection for the alternating direction method of multipliers (admm):
  quadratic problems,'' \emph{IEEE Transactions on Automatic Control}, vol.~60,
  no.~3, pp. 644--658, 2015.

\end{thebibliography}
\begin{IEEEbiography}[{\includegraphics[width=2.54cm,height=3.18cm]{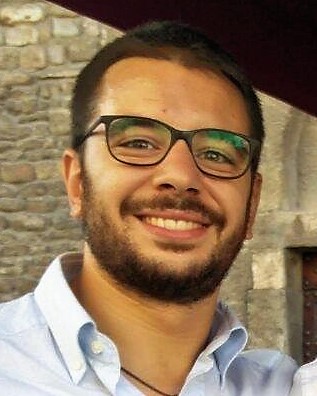}}]{Alessio La Bella}
 received the B.Sc. and the M.Sc. in Automation and Control engineering from Politecnico di Milano in 2013 and in 2015, respectively. In 2015 he received the Alta Scuola Politecnica Diploma together with the M.Sc. in Mechatronics from Politecnico di Torino. Currently, he is enrolled as a Ph.D. student in Information Technology at Politecnico di Milano. His main interests concern the design of optimization-based control algorithms for the management of energy and electrical systems, e.g. micro-grids and renewable sources.
\end{IEEEbiography}
\vspace{-10mm}
\begin{IEEEbiography}[{\includegraphics[width=2.54cm,height=3.18cm]{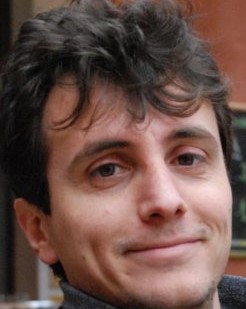}}]{Marcello Farina}
received the Laurea degree in Electronic Engineering in 2003 and the PhD degree in Information Engineering in 2007, both from the Politecnico di Milano. In 2005 he was visiting student at the Institute for Systems Theory and Automatic Control, Stuttgart, Germany. He is presently Associate Professor at Dipartimento di Elettronica, Informazione e Bioingegneria, Politecnico di Milano. His research interests include distributed, hierarchical and decentralized state estimation and control, stochastic model predictive control, and applications, e.g., mobile robots, sensor networks, and energy supply systems.
\end{IEEEbiography}
\vspace{-5mm}
\begin{IEEEbiography}[{\includegraphics[width=2.54cm,height=3.18cm]{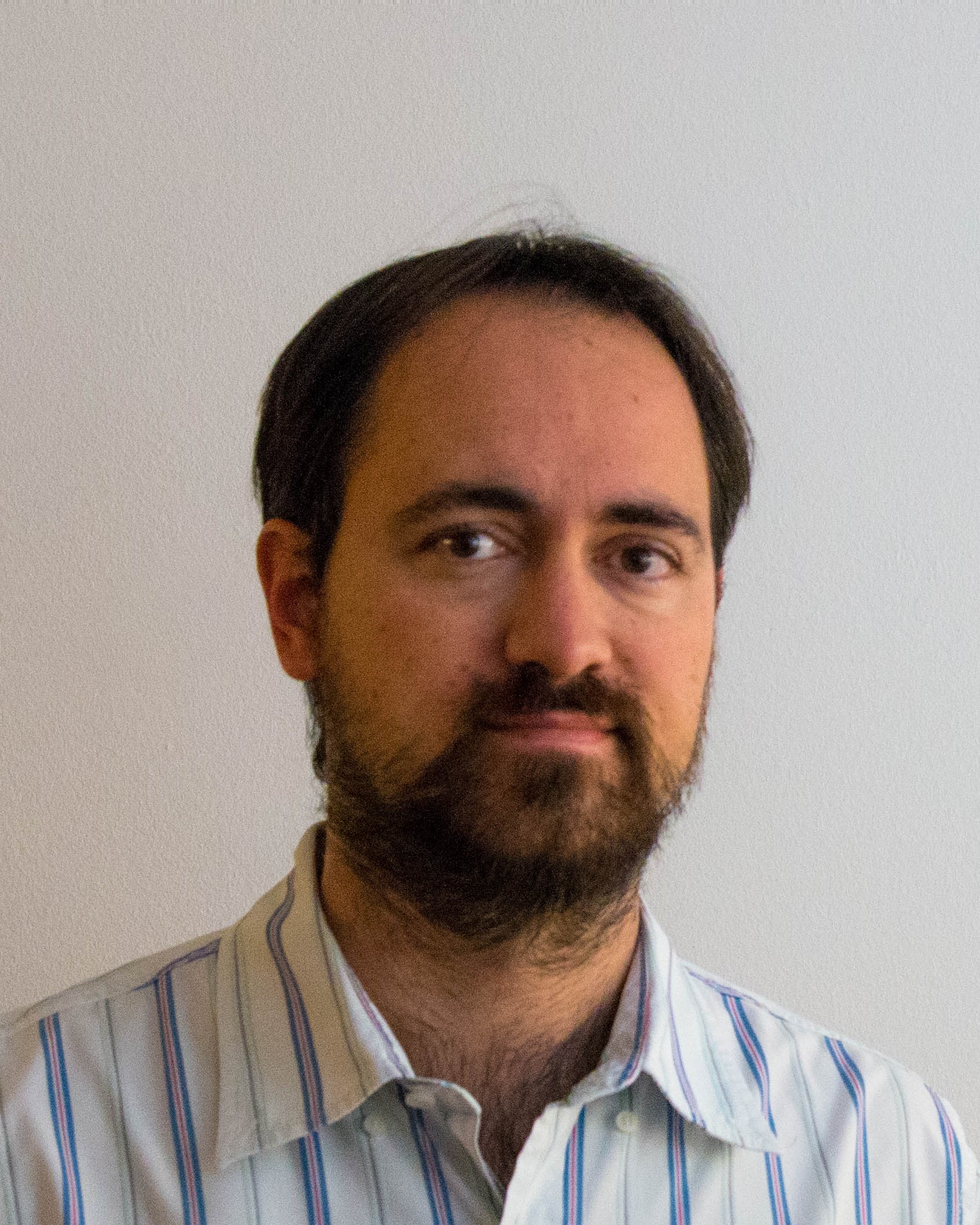}}]{Carlo Sandroni}
graduated from Politecnico di Milano with a M.Sc. in Control Engineering in 2008. Later in the same year he joined RSE in the Power Generation Technologies and Materials Department and he is currently leading the Active Networks, Automation and Control Research Group. His research topics include voltage control in distribution networks, optimization of distributed energy resources, simulation of distribution networks and generation plants, control of generators, predictive control and multi-agent systems. He has participated in several national and European projects, including DERRI, ERIGrid and ELECTRA.
\end{IEEEbiography}
\vspace{-5mm}
\begin{IEEEbiography}[{\includegraphics[width=2.54cm,height=3.18cm]{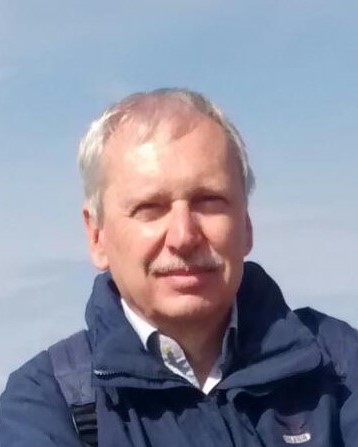}}]{Riccardo Scattolini}
is Full Professor of Automatic Control at the Politecnico di Milano, Italy. He was awarded Heaviside Premium of the Institution of  Electrical Engineers. United Kingdom and is Associate Editor of the IFAC journal Automatica. His main research interests include modeling, identification, simulation, diagnosis, and control of industrial plants and electrical systems, with emphasis on the theory and applications of Model Predictive Control and fault detection methods to large-scale and networked systems.
\end{IEEEbiography}
\end{document}